\documentclass[aap,seceqn,MSNbibl,citesort,dvips]{arximspdf}

\doi{10.1214/10-AAP752}
\volume{21}
\issue{5}
\pubyear{2011}
\firstpage{1627}
\lastpage{1662}

\makeatletter

\def\bptnote#1{}

\newtheorem{theorem}{Theorem}[section]
\newtheorem{prop}{Proposition}[section]

\newproclaim{rem}{Remark}[section]
\newproclaim{example}{Example}[section]

\newcommand{\indic}{\mathbb{I}}
\newcommand{\rme}{\mathrm{e}}

\makeatother

\begin{document}
\begin{frontmatter}

\title{Quantile clocks}
\runtitle{Quantile clocks}

\begin{aug}
\author[A]{\fnms{Lancelot F.} \snm{James}\corref{}\thanksref{t1}\ead[label=e1]{lancelot@ust.hk}} and
\author[B]{\fnms{Zhiyuan} \snm{Zhang}\ead[label=e2]{ismtzzy@gmail.com}}
\runauthor{L. F. James and Z. Zhang}
\affiliation{Hong Kong University of Science and Technology}
\address[A]{Department of Information Systems\\
Business Statistics\\
\quad and Operations Management\\
Hong Kong University\\
\quad of Science and Technology\\
Clear Water Bay, Kowloon\\
Hong Kong\\
\printead{e1}}
\address[B]{School of Statistics and Management\\
Shanghai University of Finance and Economics\\
No. 777 Guoding Road\\
Shanghai 200433\\
China\\
\printead{e2}}%adresu isvedimo komanda gale!
\end{aug}

\thankstext{t1}{Supported in part by the Grant RGC-HKUST 600907 of the HKSAR.}

% HISTORY:
\received{\smonth{3} \syear{2010}}
\revised{\smonth{8} \syear{2010}}

% ABSTRACT
%
\begin{abstract}
Quantile clocks are defined as convolutions of subordinators $L$, with
quantile functions of positive random variables. We show that quantile
clocks can be chosen to be strictly increasing and continuous and
discuss their practical modeling advantages as \textit{business activity
times} in models for asset prices. We show that the marginal
distributions of a quantile clock, at each fixed time, equate with the
marginal distribution of a single subordinator. Moreover, we show that
there are many quantile clocks where one can specify~$L$, such that
their marginal distributions have a desired law in the class of
generalized $s$-self decomposable distributions, and in particular the
class of self-decomposable distributions. The development of these
results involves elements of distribution theory for specific classes
of infinitely divisible random variables and also decompositions of a
gamma subordinator, that is of independent interest. As applications,
we construct many price models that have continuous trajectories,
exhibit volatility clustering and have marginal distributions that are
equivalent to those of quite general exponential L\'{e}vy price models.
In particular, we provide explicit details for continuous processes
whose marginals equate with the popular VG, CGMY and NIG price models.
We also show how to perfectly sample the marginal distributions of more
general classes of convoluted subordinators when $L$ is in a sub-class of
generalized gamma convolutions, which is relevant for pricing of
European style options.
\end{abstract}

% KEYWORDS
%

\pdfkeywords{60E07, 60G09, Generalized gamma convolutions,
Levy processes, perfect sampling, self-decomposable laws,
time changed price processes}

\end{frontmatter}

%s1 ###
\section{Introduction}\label{intro}
Let $Q_{R}(u)=\inf\{t\dvtx F_{R}(t)\ge u\}, 0<u<1$ denote the quantile
function of a nonnegative continuous random variable $R$ with
strictly increasing cumulative distribution function (c.d.f.) $F_{R}$, and finite
first moment $\mathbb{E}[R]$. In this paper, we introduce and
describe detailed distributional properties of a class of random
time changes $T_{R}:=(T_{R}(t), t\ge0)$, which we call
\textit{quantile clocks}. These processes are defined as,
%
%e1.1 ###
\begin{equation} \label{QQ}
T_{R}(t)=\int_{0}^{t}Q_{R}\bigl((1-s/t)_{+}\bigr)L(ds),\qquad t\ge0,
\end{equation}
where $L$ is a subordinator, and $(a)_{+}:=\max(0,a)$. While
applicable in many settings, we
follow the framework in \cite{Bender} and discuss the modeling
advantages of quantile clocks as \textit{business activity times} in
time changed models for asset prices.

The quantile clocks, may be written as special cases of
\textit{convoluted subordinators}, which are processes described in
Bender and Marquardt \cite{Bender}. That is to say, processes
defined as $T(t)=\int_{0}^{t}k(t,s)L(ds), t\ge0$, for $k(t,s)$ a
known kernel. The authors \cite{Bender}, Proposition 1, provide
mild conditions on $k(t,s)$ and $L$ such that a process
$T:=(T(t),t\ge0)$ has almost surely strictly continuous and
increasing sample paths. In terms of applications, \cite{Bender}
argue that one can use $W(T(t))$, where $W$ is a brownian motion
with drift, as time changed models for the log price of assets
that possess continuous trajectories, where $T(t)$ is now
interpreted as \textit{business activity time}. Furthermore, such
models can correct deficiencies in Black--Scholes type price
models. In particular, it is known that (i) the log returns of
asset prices have nonnormal distributions, and often exhibit
semi-heavy or heavier tail behavior, (ii) the volatility or
variance is dependent on time, (iii) asset prices exhibit volatility
\textit{clustering} or \textit{persistence.} Reference \cite{Bender} also describe
a general formula for
European style option prices that depend on the marginal
distribution of $T(\tau)$ for some fixed time to maturity $\tau>0$.
For other applications of processes representable as convoluted
subordinators, see, for instance,
\cite{James05,PeccatiPrunster,PeccatiTaqqu08,WolpertIckstadt} and references therein.

In the literature, exponential L\'{e}vy price models, defined as
${\rme} ^{-\chi(t)}$ for a L\'{e}vy process $\chi$ on $\mathbb{R}$,
have been quite successful in terms of their ability to capture
some of the stylistic features of asset prices (i) and (ii) listed
above. In addition, there are many choices of $\chi$ where one can easily
calibrate pricing models to the options market, capturing
volatility smiles and skews, via Monte Carlo methods or perhaps
more generally by the fast Fourier transform (FFT) methods
outlined in Carr and Madan \cite{CarrFFT}. Many L\'{e}vy
processes $\chi$, can be expressed as $W(\zeta(t))$ for
some subordinator $\zeta$. However, the precise $\zeta$
that is associated with a $\chi$ is not always known explicitly, and
$\chi$ is often modeled via its L\'{e}vy density. Arguably, the most
popular models of this type include the variance gamma process
(VG) by \cite{MadanVG}, where $\zeta$ is a gamma subordinator, the
Carr--Geman--Madan--Yor (CGMY) process \cite{CGMY}, where $\zeta$ has
only recently been identified by Madan and
Yor \cite{MadanYorCGMY}, and the normal inverse Gaussian (NIG)
process \cite{BN98}, where $\zeta$ is an inverse Gaussian process.
The popularity of these models is due to their relative simplicity
and distributional flexibility.

However,
since $\chi$ has independent increments, exponential L\'{e}vy
processes are unable to capture effects due to volatility
clustering. One approach discussed in \cite{CarrLevy}, which is
related to the price models in \cite{BNS2001}, is to
further time change $\chi$ by a stochastic volatility process of
the form $T(t)=\int_{0}^{t}v(s)\,ds$, where $v(s)$ represents the
instantaneous volatility either following a mean reverting
Cox--Ingersoll--Ross (CIR) process or a non-Gaussian
Ornstein--Uhlenbeck (BNS--OU) model of Barndorff-Nielsen and
Shephard \cite{BNS2001}, specified by the dynamics
%
%e1.2 ###
\begin{equation}
\label{BNSOUSDE}
\mathrm{d}v(t)=-\lambda v(t)\,dt+\vartheta(\lambda \,dt),
\end{equation}
where $\vartheta$ is a subordinator we shall call an OU--BDLP. The
BNS--OU model $v(t)$, possesses jumps, and has a stationary
distribution with laws equating to the class of laws of
self-decomposable random variables that remarkably one can choose
based on a prescribed choice of $\vartheta$. This latter fact is
important for our exposition. In contrast, the integrated
volatility $T(t)=\int_{0}^{t}v(s)\,ds$ is continuous and has
nontypical marginals laws (obviously depending on $t$). In the case
of the CIR process, $v(t)$ is a diffusion having a transition
density following a noncentral chi-squared distribution. We will
not consider models of this type.

The authors \cite{Bender} demonstrate that their approach,
involving convoluted subordinators, can be viewed as viable
variations of the idea in \cite{CarrLevy}. Furthermore, their work
essentially contains the popular model of \cite{BNS2001}. However,
as noted by the authors, there are practical issues arising for
instance in the pricing of options that relate to the marginal
distributions of $T(\tau)$ at maturity times $\tau$. While
$T(\tau)$ are infinitely divisible, their marginal distributions
and also characteristic functions depend on $k(\tau,s)$ and
$\tau$, in a nontrivial way. Hence, leading in general to
nonfamiliar distributions for $T(\tau)$. Related to this point,
is a classical problem where in general it is not clear how to
\textit{exactly} sample infinitely divisible random variables, even
in the case of a subordinators $\zeta(\tau)$ for each fixed
$\tau$. Some notable exceptions for $\zeta$ are gamma, positive
stable and inverse Gaussian processes whose marginals for each
fixed time point are gamma, positive stable and inverse Gaussian
random variables, and hence are easily sampled. More generally, one
can resort to sampling methods based on truncation of infinite
series representations, but these do not yield exact samples and
it is not always clear how to control the level of accuracy.
Reference \cite{Bender} do point out that if $T(\tau)$ has an analytically
tractable characteristic function then one can apply the popular
fast Fourier transform (FFT) techniques in Carr and
Madan \cite{CarrFFT} to obtain explicit option prices. They also
provide supporting results for some choices of $k$ and $L$.
However, due to the generality of $k(t,s)$, it appears difficult
to apply the (FFT) for option price formula depending on general
$T(\tau)$. These points do not reflect a deficiency in the
approach of \cite{Bender} but rather that the class of convoluted
subordinators is quite general. The task then becomes how to
choose kernels $k$ and subordinators $L$ that are convenient in
terms of implementation as well as having general modeling
flexibility.

%s1.1 ###
\subsection{Contributions and outline}\label{sec1.1}
In Section \ref{sec3}, we show that quantile clocks, which arise by setting
$k(t,s)=Q_{R}((1-s/t)_{+})$, can be chosen to have continuous
trajectories and in general have marginal distributions that, for
each fixed time $t$, equate to a single subordinator $\zeta$. That
is $T_{R}(t)\stackrel{d}{=}\zeta(t)$ for fixed $t$. We also highlight
a very tractable example related to \cite{JamesLamperti,JRY}. Of
course, in general, the law $\zeta$, that is, the marginal laws of
$T_{R}$, depends on $(R,L)$, and hence the deterministic quantile
function $Q_{R}$. However, in Section \ref{sec4}, we show that there are
many quantile clocks whose marginal distributions can be chosen
such that they do not obviously depend on $Q_{R}$. In fact for
these given $Q_{R}$, we show that one can choose $L$, and a random
variable $Y$, such that the marginals of $T_{R}$ have specific
laws in the Jurek's \cite{IksanovJurek,Jurek85,Jurek88,Jurek89}
$\mathcal{U}_{\delta}$ class of \textit{generalized}
$s$-selfdecomposbale laws, for $\delta>0$. These classes
contain the important class $\mathcal{L}$ of self-decomposable
distributions on $\mathbb{R}_{+}$. See \cite{CarrSD} for the
relevance of self-decomposable L\'{e}vy processes in financial
modeling.

This ability to choose specific (familiar) marginal laws for price
processes, while allowing for quite varied path properties induced
by different quantiles~$Q_{R}$, gives modelers a great deal of
flexibility. It is also reminiscent of how one might choose\vspace*{-1pt} a
BNS--OU model $v$ to have a specific stationary distribution, that is,
$v(0)\stackrel{d}{=}v(t)$ for all fixed $t$, in $\mathcal{L}$ based
on the OU--BDLP $\vartheta$ appearing in  (\ref{BNSOUSDE}).
However, recall that $v$ has jumps and the law of $T_{R}$
obviously must depend on $t$. The precise methods we use to
establish these results, and identify $Q_{R}$, $Y$ and $L$, are
given in Sections \ref{sec4}, \ref{sec5} and \ref{sec6}, and\vspace*{-1pt} should also be of general
interest to experts in L\'{e}vy processes. In Section \ref{sec7}, we exploit
the fact that $T_{R}(t)\stackrel{d}{=}\zeta(t)$ for each $t$, and we
show that compositions (or time changes) involving quantiles
clocks behave marginally like subordinators. In Section \ref{sec8}, we
show that as consequences of our results, that we are able to
identify price processes whose marginal behavior coincides with
those of exponential L\'{e}vy price processes. In particular, we
identify explicitly many processes whose marginal distributions
are equivalent to VG, CGMY and NIG price processes, but whose
trajectories are continuous and otherwise quite varied, and
additionally exhibit volatility clustering. We also identify
models possessing jumps that otherwise have the properties
mentioned above. In Section \ref{sec9}, we show how one can use our results
for quantile clocks to specify laws for the convoluted
subordinator referred as a short memory kernel in \cite{Bender}.

While quantile clocks are our main focus, in Section \ref{sec2}, we also
describe results that apply to the practical implementation of log
price models $W(T(t))$, considered in \cite{Bender}, where $T(t)$
is based on a general kernel $k(t,s)$. In particlar, if $L$ is
chosen to have laws in the class of generalized gamma convolutions
with finite Thorin measure,
see \cite{BondBook,JamesBernoulli,JRY}, call this class
$\mathcal{G}_{+}$, then the random variable $T(\tau)$ can be
exactly sampled in many instances. This is based on a very recent
work of Devroye and James \cite{DevroyeJames} where a double
coupling from the past (Double CFTP) perfect sampling routine is
devised, and also results described in
James \cite{JamesBernoulli}. Furthermore, for this choice of $L$,
by using a deterministic time-change we obtain a simplified
version of the option price formulae given in \cite{Bender}.

%s1.2 ###
\subsection{Preliminaries}\label{sec1.2}
We now present some concepts and notation we shall use throughout.
First, for fixed positive numbers $(a,b)$, let $\gamma_{a}$
denote a $\operatorname{gamma}(a)$ random variable with shape parameter
$a$ and scale $1$, let $\beta(a,b)$ denote a beta $(a,b)$
random variable. Furthermore,\vspace*{-1pt} $U$ will always denote a
Uniform$[0,1]$ variable, and recall that for any
$\delta>0$, $U^{1/\delta}\stackrel{d}{=}\beta_{\delta,1}$. $\xi_{p}$
is a Bernoulli random variable with success probability $p$. In
addition for a generic random variable $Y$, $Y'$ will denote a
variable equivalent in distribution but otherwise independent. $(N(s)\dvtx s>0)$
will denote a homogeneous Poisson process with intensity $\mathbb{E}[N(s)]=s$.
For
a (nonrandom) function $g(x)$, $g'(x)$ and $g''(x)$ with denote its
first and
second derivatives.

Formally, recall that a subordinator $\zeta=(\zeta(t);t>0)$, is an
increasing process with right continuous paths and stationary
independent increments, whose law is specified by its Laplace
transform for some $\omega>0$
%
%e1.3 ###
\begin{equation}
\label{LaplaceID}
\mathbb{E}\bigl[{\rme}^{-\omega\zeta(t)}\bigr]={\rme}^{-t\psi_{\zeta
}(\omega)},
\end{equation}
where for some $c\ge0$, $\omega>0$
\[
\psi_{\zeta}(\omega)=c\omega+\int_{0}^{\infty}(1-{\rme}^{-\omega s})\Lambda_{\zeta}(ds)
\]
is finite and is called the \textit{Laplace exponent} of $\zeta$,
$\Lambda_{\zeta}(ds)$ is its L\'{e}vy measure,
$\rho_{\zeta}(s)=\Lambda_{\zeta}(ds)/ds$ is the L\'{e}vy density. We
will work with the case where $c=0$.

It follows that the laws of $\zeta$ can be specified by any of
these quantities. $\zeta$ is said to be of \textit{infinite
activity} if $\Lambda_{\zeta}(\infty)=\infty$ and otherwise of
\textit{finite activity}. In the latter case, $\zeta$ corresponds to
compound Poisson process whose jumps have a common probability
density/mass function proportional to $\rho_{\zeta}$. Throughout we
shall reserve the notation $\zeta, L,Z $ for generic
subordinators, and corresponding random variables, and the
notation $\vartheta$ for the OU--BDLP. As is well known, for each
fixed $t$, $\zeta(t)$ is a random variable in the class
$\mathcal{J}$ of infinitely divisible random variables (taking
values in $\mathbb{R}_{+}$). We now describe the characteristics
of some important subclasses of $\mathcal{J}$, say $\mathcal{L},
\mathcal{B}, \mathcal{G}$ and $\mathcal{G}_{+}$, satisfying $
\mathcal{G}_{+}\subset\mathcal{G} \subset\mathcal{L}\subset
\mathcal{J}$ and $ \mathcal{G}_{+}\subset\mathcal{G} \subset
\mathcal{B} \subset\mathcal{J}$.

We say that a random variable $\zeta(1)$ is in the class
$\mathcal{L}$ of self-decomposable variables if
$\rho_{\zeta}(s)=s^{-1}h(s)$, with $h$ decreasing. We also note
that from Jurek and Vervaat \cite{JV} that, with respect to the
OU process in (\ref{BNSOUSDE}), there is the relationship
\[
\zeta(1)\stackrel{d}{=}v(0)\stackrel{d}{=}v(t)\stackrel{d}{=}\int
_{0}^{\infty}{\rme}^{-s}\vartheta(ds).
\]
Note that we will say that $\zeta$ is a subordinator in
$\mathcal{L}$ to mean that it is a subordinator whose L\'{e}vy
density corresponds to that of a variable in $\mathcal{L}$, of
course $\zeta(t)$ is in $\mathcal{L}$ for each fixed $t$. Similar
statements will apply for other classes. We say that $\zeta(1)$ is
a variable in Bondesson's \cite{BondBook}, Section 9, $\mathcal{B}$
class, or the class of generalized convolutions of mixtures of
exponential distributions (GCMED), if the L\'{e}vy density is
completely monotone, that is,
$\rho_{\zeta}(s)=\int_{0}^{\infty}{\rme}^{-sy}\mu(dy)$, for
some nonnegative measure $\mu$.

We now describe the classes $\mathcal{G}$ and $\mathcal{G}_{+}$.
$\zeta(1)$ is a variable that is a generalized gamma convolutions
(GGC), see \cite{BondBook}, if it is in the class $\mathcal{G}$,
characterized by
\[
\rho_{\zeta}(s)=s^{-1}\int_{0}^{\infty}{\rme}^{-sy}\nu(dy)\quad
\mbox{and}\quad\psi_{\zeta}(\omega)=\int_{0}^{\infty}\log(1+\omega/y)\nu(dy)
\]
for some sigma-finite measure $\nu$, formally known as a Thorin
measure. We say that $\zeta$ is a $\operatorname{GGC}(\nu)$ subordinator.

A $\zeta(1)$ variable is in the class $\mathcal{G}_{+}$, if it
satisfies
%
%e1.4 ###
\begin{equation}
\label{gcon}
\rho_{\zeta}(s)=\theta s^{-1}\mathbb{E}[{\rme}^{-s/R}]
\quad\mbox{and}\quad\psi_{\zeta}(\omega)=\theta\mathbb{E}[\log(1+\omega R)]
\end{equation}
for some $\theta>0$ and some random variable $R$ satisfying
$\mathbb{E}[\log(1+\omega R)]<\infty$.
In this case, we say $\zeta(1)$ is a $\operatorname{GGC}(\theta,R)$
variable. Moreover, $\zeta(t)$ is a $\operatorname{GGC}(\theta t,R)$ variable for
each fixed $t$, and $\zeta$ is referred to as a $\operatorname{GGC}(\theta,R)$
subordinator.

We now highlight some important properties of $\operatorname{GGC}(\theta,R)$
random variables, and subordinators that for instance allow them
to be exactly sampled by the methods in \cite{DevroyeJames}. These
facts can found in \cite{JamesBernoulli} as well as
\cite{JRY,JamesLamperti}, and depend heavily on the results for Dirichlet
means in \cite{CifarelliRegazzini}. Letting\vspace*{-1pt} $Z_{\theta}$ denote a
$\operatorname{GGC}(\theta,R)$ subordinator, it follows that
$Z_{\theta}(t)\stackrel{d}{=}Z_{\theta t}(1)$. Importantly, there is
the representation, for any $\kappa\ge\theta>0$,
$Z_{\theta}(1)\stackrel{d}{=}\gamma_{\theta}M_{\theta}=\gamma
_{\kappa}\tilde{M}_{\kappa}$
where
%
%e1.5 ###
\begin{equation}
\label{recursionM}
M_{\theta}\stackrel{d}{=}\beta_{\theta,1}M_{\theta}+(1-\beta
_{\theta,1})R
\end{equation}
and
%
%e1.6 ###
\begin{equation}
\label{recursionM2}
\tilde{M}_{\kappa}\stackrel{d}{=}\beta_{\kappa,1}\tilde
{M}_{\kappa}+(1-\beta_{\kappa,1})R\xi_{p}
\end{equation}
for $p=\theta/\kappa$. That is, a $\operatorname{GGC}(\theta, R)$ random variable
variable is a $\operatorname{GGC}(\kappa,R\xi_{p})$ variable. In particular, if
$0<\theta t\leq1$, then $Z_{\theta
t}(1)=Z_{\theta}(t)\stackrel{d}{=}\gamma_{1}\tilde{M}_{1}$, where
%
%e1.7 ###
\begin{equation}
\label{Uni}
\tilde{M}_{1}\stackrel{d}{=}U\tilde{M}_{1}+(1-U)R\xi_{p}
\end{equation}
for $p=\theta t$, and it follows from the work of Cifarelli and
Regazzini \cite{CifarelliRegazzini} that $\tilde{M}_{1}$ has
density of the form
\[
\frac{x^{p-1}}{\pi}\sin(\pi F_{R\xi_{p}}(x))
{\rme}^{-p\Psi_{R}(x)}\qquad\mbox{for }x>0
\]
with
\[
\Psi_{R}(x)=\mathbb{E}\bigl[{\log}|x-R|\indic_{(R\neq x)}\bigr].
\]
Thus, as pointed out in \cite{JamesBernoulli}, if one one can
evaluate $\Psi_{R}(x)$ in a suitable fashion, then one can exactly
sample any variable $Z_{\theta}(t)$ for every fixed $0<t\leq
1/\theta$, by, for instance, rejection sampling. Since any number
$s>0$, can be set to $s=nt$, for some integer $n$ and $0<t\leq
1/\theta$, it follows that
$Z_{\theta}(s)\stackrel{d}{=}Z_{\theta}(tn)$ can be exactly sampled
by at most exactly sampling $n$ copies of the random variable
$Z_{\theta}(t)\stackrel{d}{=}\gamma_{1}\tilde{M}_{1}$. We note that
in general $M_{\theta}$ for $\theta>0$, does not have a simple
expression for its density. So the exact sampling method suggested
above relies solely on the ability to sample the variable in
$\tilde{M}_{1}$ in (\ref{Uni}), for each $p$. This is possible
provided that $\Psi_{R}(x)$ is analytically tractable. However,
since $R$ can be quite arbitrary this will not always be true.
Fortunately, there is the recent Double CFTP perfect sampling
method by \cite{DevroyeJames} that can be used to exactly sample
any of the variables satisfying (\ref{recursionM}), (\ref
{recursionM2}) or (\ref{Uni}).
This procedure applies provided
that $R$ is a bounded variable and one has a method to sample $R$,
but otherwise does not require any potentially complicated
calculations. Hence, any $\operatorname{GGC}(\theta, R)$ variable, with $R$
bounded, can be exactly sampled by drawing an independent gamma
variable and applying the Double CFTP. Details may be found
in \cite{DevroyeJames}, however we shall sketch out the details
for a subclass of the variables $T(\tau)$ in the next section.
\begin{rem}
Letting $Q_{R}(u)$ denote a quantile function of $R$, variables
$M_{\theta}$, satisfying (\ref{recursionM}), are called Dirichlet
means since they can always be represented as
\[
M_{\theta}\stackrel{d}{=}\int_{0}^{1}Q_{R}(u)D_{0,\theta
}(du|F_{U})\stackrel{d}{=}\int_{0}^{\infty}yD_{\theta}(dy|F_{R}),
\]
where
\[
D_{\theta}(y|F_{R})\stackrel{d}{=}\sum_{k=1}^{\infty}P_{k}\indic
_{(R_{k}\leq
y)}
=\sum_{k=1}^{\infty}V_{k}\prod_{j=1}^{k-1}(1-V_{j})\indic
_{(R_{k}\leq
y)}
\]
is a Dirichlet process with $(P_{k})$ a sequence of
probabilities having a Poisson Dirichlet law with parameter
$\theta$, see \cite{Ferg73,PY97}. That is, for each $k$,
$\gamma_{\theta}P_{k}\stackrel{d}{=}J_{k}$, where $(J_{k})$ are the
ranked jumps of a $\operatorname{gamma}(\theta)$ subordinator. $(V_{k})$ are i.i.d.
Beta$(1,\theta)$ random variables, and $(R_{k})$ are i.i.d. $F_{R}$.
See \cite{JRY,LijoiMean} for more details.
\end{rem}
\begin{rem}\label{stableremark}For $0<\alpha<1$, positive stable subordinators
$S_{\alpha}(t)$, where $S_{\alpha}(1):=S_{\alpha}$, with
$\psi_{S_{\alpha}}(\omega)=\omega^{\alpha}$, and corresponding
processes $\widehat{S}_{\alpha}(t)$, with\break
$\psi_{\widehat{S}_{\alpha}}(\omega)=(1+\omega)^{\alpha}-1$, as well
as\vspace*{-1pt} their scaled variations, are in $\mathcal{G}$ but not $\mathcal
{G}_{+}$. Naturally a $\operatorname{gamma}(\theta)$ subordinator, say
$(\gamma_{\theta}(t);t\ge0)$, is in $\mathcal{G}_{+}$. However,
$S_{\alpha}, \widehat{S}_{\alpha}$ and $\gamma_{\theta}$,
constitute a family of (generalized gamma) subordinators
with L\'{e}vy density
\[
Cs^{-\alpha-1}e^{-bs}
\]
for $0\leq\alpha<1$ and $b\ge0$, see \cite{PY97}, Proposition 21.
Additionally, heavy tailed
variables such as Linnik variables of the form
$S_{\alpha}(\gamma_{\theta
t})\stackrel{d}{=}\gamma^{1/\alpha}_{\theta t}S_{\alpha}$ are in
$\mathcal{G}_{+}$. As well as their exponentially tilted
counterparts $\widehat{S}_{\alpha}(\gamma_{\theta t}p)$, for some $0<p<1$.
See \cite{JamesLamperti}.
\end{rem}

%s2 ###
\section{Convoluted subordinators}\label{sec2}
We now give the formal specifications for convoluted subordinators
as defined in Bender and Marquardt \cite{Bender}. Throughout the
rest of the paper, let $L$ denote an infinite activity
subdordinator. That is the L\'{e}vy measure,
$\Lambda_{L}(\infty)=\infty$. In order that the convoluted
subordinator
%
%e2.1 ###
\begin{equation}
T(t)=\int_0^tk(t,s)\,{d}L(s)
\end{equation}
has strictly continuous and increasing trajectories, $k$ is chosen
to satisfy the following regularity conditions:
\begin{longlist}[(a)]
\item[(a)] for fixed $t\in[0,\infty)$, the mapping $s\mapsto k(t,s)$
is integrable,
\item[(b)] for fixed $s\in[0,\infty)$, the mapping $t\mapsto k(t,s)$
is continuous and increasing and there is an $\varepsilon>0$ such that
$t\mapsto k(t,s)$ is strictly increasing on $[s,s+\varepsilon]$,
\item[(c)] $k(t,s)=0$ whenever $s>t\geq0$.
\end{longlist}

The authors also derive a weighted Black--Scholes pricing formula
for European style options as follows. Let
%
%e2.2 ###
\begin{equation}
\label{drift} \widehat{W}_{\mu}(t)=W(t)+\mu t
\end{equation}
denote a standard Brownian motion with drift parameter $\mu$, that
is $W(t)$ is a standard Brownian motion. Recall that for geometric
Brownian motion the price process under the risk neutral measure is
given by
%
%e2.3 ###
\begin{equation}
\label{priceBS} S(t)=S(0)\exp\{rt+\widehat{W}_{-1/2}(\sigma^{2}t)\},
\end{equation}
where $\widehat{W}_{-1/2}$ is defined by
(\ref{drift}) with $\mu=-1/2$.

Setting $S_{\tau}=S(\tau)$, the quantity $(S_{
\tau}-K)_{+}$ is the the payoff function of a European call option
with strike $K>0$ and maturity $\tau$, and $r>0$ is the risk-free
interest rate. Then the Black--Scholes formula for the price at time
$0$, say $\mathrm{B}(\sigma, K,\tau)$, is given by
%
%e2.4 ###
\begin{equation}
\label{BSprice}\quad
\mathrm{B}(\sigma, K,\tau)=
{\rme}^{-r\tau}\mathbb{E}[(S_{\tau}-K)_{+}]=S_{0}\Phi(d_{1}(
\sigma))-K{\rme}^{-r\tau}\Phi(d_{2}(
\sigma)),
\end{equation}
where $\Phi(x)$ is the standard normal distribution function
\[
d_{1}(\sigma)=\frac{\log(S_{0}/K)+(r+\sigma^{2}/2)\tau}{\sigma
\sqrt{\tau}}\quad\mbox{and}\quad d_{2}(\sigma)=d_{1}(\sigma)-\sigma\sqrt
{\tau}.
\]
See Schoutens \cite{Schoutens} for this notation.
In \cite{Bender}, a
price model under risk neutral dynamics is specified as
%
%e2.5 ###
\begin{equation}
\label{priceCS} \tilde{S}(t)=S_{0}\exp\{rt+\widehat{W}_{-1/2}(
\sigma^{2}T(t))\},
\end{equation}
where now $\tilde{S}$ is the asset price, and $T(t)$ is a convoluted
subordinator. They obtain the following pricing formula.
\begin{theorem}[(Bender and
Marquardt \cite{Bender}, Theorem 4)]\label{tina} For the price model (\ref
{priceCS}), with $\tilde{S}_{\tau}=\tilde{S}(\tau)$.
Let $(\tilde{S}_{\tau}-K)_{+}$ be the payoff\vspace*{1pt} function of a
European call option with strike $K\in\mathbb{R}_{+}$ and
maturity $\tau$. Then the initial fair price of $(\tilde{S}_{\tau
}-K)_{+}$, is given by
%
%e2.6 ###
\begin{equation}
\label{Benderprice}
{\rme}^{-r\tau}\mathbb{E}[(\tilde{S}_{\tau}-K)_{+}]=\mathbb
{E}\bigl[\mathrm{B}\bigl(\sigma\sqrt{T(\tau)/\tau}, K,\tau\bigr)\bigr],
\end{equation}
where for positive $y$, $\mathrm{B}(y, K,\tau)$ is the Black--Scholes
price given in (\ref{BSprice}) with $y$ in place of
$\sigma$. Furthermore, $S_{0}$ is considered fixed.
\end{theorem}

As noted in \cite{Bender}, and discussed in the \hyperref[intro]{Introduction}, the
problem with the above result is that it is in general difficult
to handle the exact law of $T(\tau)$. However, the authors do
point out that if $T(\tau)$ possesses an analytically tractable
characteristic function then it is possible to use fast Fourier
transform (FFT) methods. They give some special examples where this
might be possible, but in general this is not straightforward.
This is clear since the Laplace exponent of $T(\tau)$ can be
expressed as
%
%e2.7 ###
\begin{equation} \label{LCS}
\psi_{T(\tau)}(\omega)=\tau\mathbb{E}[\psi_{L}(\omega
k(\tau,U\tau))],
\end{equation}
where $k$ and $L$ are quite general and the expression otherwise
depends on $\tau$ in a nontrivial way.

We believe Theorem \ref{tina} does have quite a bit of utility
provided that one can have more control over the choice of
marginal laws exhibited by $T(\tau)$, for each fixed $\tau$. Next,
we show that by choosing $L$ to be in $\mathcal{G}_{+}$ one can
(in a practical sense) use Theorem \ref{tina} for many kernels
$k$. Even those that do not admit nice characteristic functions.
\begin{rem}
We note that the martingale argument used in \cite{Bender} is
different than that used for standard time changed models.
The filtration used by \cite{Bender} preserves the martingale property
for a larger class of models including, of course, time changes by a
simple subordinator. However, the usual filtrations used for simple
subordinators may not preserve the martingale property for all
convoluted subordinators.
\end{rem}

%s2.1 ###
\subsection{A general result for $L$ in $\mathcal{G}_{+}$}\label{sec2.1}

As we just mentioned, we now look at the choice where $L$ is a
$\operatorname{GGC}(\theta,Y)$ subordinator where $Y$ is some random variable. In
terms of modeling for general $T(t)$, we can view $\theta$ as a
time parameter that can be manipulated for practical convenience.
This is partly because the variable $Y$ can have unknown
parameters that can be used for calibration.
\begin{theorem}\label{genGGC}Let $\widehat{W}_{\mu}(t)$ denote a
Brownian motion with drift as defined
in (\ref{drift}). Let
%
%e2.8 ###
\begin{equation}
\label{GGCclock}
T(t)=\int_{0}^{t}k(t,y)L_{\theta}(dy)
\end{equation}
denote a
convoluted subordinator where $L_{\theta}$ is a $\operatorname{GGC}(\theta, Y)$
subordinator. For each fixed $t$, define a random variable $
R_{t}\stackrel{d}{=}k(t,Ut)$. Then the process $
(\widehat{W}_{\mu}(T(t))\dvtx t\ge0)$ is almost surely continuous and
has the following
distributional properties:
\begin{longlist}[(iii)]
\item[(i)] for each fixed $t$, $T(t)$ is a $\operatorname{GGC}(\theta t, R_{t}Y)$
random variable satisfying $T(t)\stackrel{d}{=}\gamma_{\theta
t}M_{\theta t}$, where
%
%e2.9 ###
\begin{equation}
\label{repM}
M_{\theta t}=\beta_{\theta t,1}M_{\theta t}+(1-\beta_{\theta
t,1})R_{t}Y,
\end{equation}
\item[(ii)]if $0<\theta t=p\leq1$, then $T(t)\stackrel{d}{=}\gamma
_{1}M_{1,t}$,
where
%
%e2.10 ###
\begin{equation}
\label{shortM}
M_{1,t}\stackrel{d}{=}UM_{1,t}+(1-U)R_{t}Y\xi_{p},
\end{equation}
\item[(iii)] the density of the ${M}_{1,t}$
is given by
\[
\frac{x^{p-1}}{\pi}\sin(\pi\tilde{F}_{t}(x))
{\rme}^{-p\tilde{\Psi}_{t}(x)}\qquad\mbox{for }x>0,
\]
where
$\tilde{F}_{t}=F_{R_{t}Y\xi_{p}}$ is the c.d.f. of the variable
$R_{t}Y\xi_{p}$ and
\[
\tilde{\Psi}_{t}(x)=\mathbb{E}\bigl[{\log}|x-R_{t}Y|\indic_{(R_{t}Y\neq
x)}\bigr].
\]
\end{longlist}
\end{theorem}

\begin{pf}
Following (\ref{LCS}) and (\ref{gcon}), it is easy to see that the
Laplace exponent of $T(t)$ is given by
\[
\psi_{T(t)}(\omega)=t \mathbb{E}[\psi_{L_{\theta}}(\omega
k(t,Ut))]=\theta t \mathbb{E}[\log(1+\omega R_{t}Y)].
\]
The results (i), (ii) and (iii) then follow from the material we
discussed at the end of Section \ref{sec1.2}.
\end{pf}

We now state a result for European style options, which is immediate
from Theorems \ref{tina} and \ref{genGGC}.
\begin{prop}For the price\vspace*{1pt}
model (\ref{priceCS}), let $T(t)$ be the process specified
by (\ref{GGCclock}) and otherwise consider the setup in
Theorem \ref{tina}. Let $(\tilde{S}_{\tau}-K)_{+}$ be the payoff
function of a European\vspace*{1pt} call option with strike $K\in
\mathbb{R}_{+}$ and maturity $\tau$. Then the initial fair price
of $(\tilde{S}_{\tau}-K)_{+}$, is now given by
%
%e2.11 ###
\begin{equation}
\label{Gammaprice}
{\rme}^{-r\tau}\mathbb{E}[(\tilde{S}_{\tau}-K)_{+}]=\mathbb{E}\bigl[\mathrm
{B}\bigl(\sigma\sqrt{\gamma_{\theta\tau}M_{\theta,\tau}/\tau},
K,\tau\bigr)\bigr].
\end{equation}
\end{prop}

The pricing formula in (\ref{Gammaprice}) can be expressed in terms of
a (VG) process with random scale $M_{\theta\tau}$ specified by
(\ref{repM}). If for $R_{\tau}=k(\tau,U\tau)$, $R_{\tau}Y$ is
bounded then one can obtain perfect samples of the distribution of
$T(\tau)$, via \cite{DevroyeJames}. For certain~$k$, that are not
necessarily bounded, one can use the density formula in (iii) of
Theorem \ref{genGGC}.

The next result introduces a nonrandom time change that leads to a
significant reduction in complexity.
First, define for $m>0$,
\[
\phi_{m}(\mu)=\sqrt{2+\mu^{2}m^{2}}/m+\mu\quad\mbox{and}\quad b_{m}(\mu
)=\frac{1}{m\sqrt{2+\mu^{2}m^{2}}}
\]
and
$c_{m}(\mu)=b_{m}(\mu)/\phi_{m}(\mu)$.
\begin{theorem} For the convoluted subordinator in Theorem \ref
{genGGC}, the time changed process
$(\tilde{X}_{\theta}(s)\dvtx s\ge0):=(\widehat{W}_{\mu}(T((1-
{\rme}^{-s})/\theta))\dvtx s\ge0)$
satisfies for each fixed $s>0$,
\[
\tilde{X}_{\theta}(s)\stackrel{d}{=}\widehat{W}_{\mu}({\gamma
_{1}M_{1,s^{*}}}),
\]
where $M_{1,s^{*}}$ satisfies (\ref{shortM}) for
$t=s^{*}=(1-{\rme}^{-s})/\theta$ and $p=1-{\rme}^{-s}$.
\begin{longlist}[(ii)]
\item[(i)]
Furthermore, for each fixed $s$, $\tilde{X}_{\theta}(s)$ given
$M_{1,s*}=m^{2}$, for $m>0$, follows a double exponential distribution,
with density
\[
f_{\tilde{X}_{\theta}(s)}(z|m)=\cases{%
b_{m}(\mu){\rme}^{z\phi_{m}(\mu)}, &\quad $z\leq0$, \cr
b_{m}(\mu){\rme}^{-z\phi_{m}(-\mu)}, &\quad $z>0$,}
\]
and distribution function
\[
F_{\mu}(z|m)=\cases{
c_{m}(\mu){\rme}^{z[\phi_{m}(\mu)]}, &\quad $z\leq0$,\cr
c_{m}(\mu)+c_{m}(-\mu)\bigl(1-{\rme}^{-z[\phi_{m}(-\mu)]}\bigr), &\quad
$z>0$.}
\]

\item[(ii)]Hence, if the price process in (\ref{Gammaprice}) is based
on substituting $T(t)$ with the time
time changed clock, $T((1-{\rme}^{-s})/\theta)$ for $s>0$, then the fair price is given by
%
%e2.12 ###
\begin{equation}
\label{Timeprice}
{\rme}^{-r\tau}\mathbb{E}[(\tilde{S}_{\tau}-K)_{+}]=\mathbb{E}[\mathrm
{DE}(\sigma^{2}M_{1,\tau^{*}},
K,\tau)],
\end{equation}
where, $\tau^{*}=(1-{\rme}^{-\tau})/\theta$, and for $z=\log
(S_{0}/K)+r\tau$,
\[
\mathrm{DE}(y^{2},
K,\tau)=S_{0}F_{-1/2}(z|y)-{\rme}^{-r\tau}KF_{1/2}(z|y).
\]
\end{longlist}
\end{theorem}
\begin{pf} The result follows from the fact that, again,
\[
\psi_{T(t)}(\omega)=t \mathbb{E}[\psi_{L_{\theta}}(\omega
k(t,Ut))]=\theta t \mathbb{E}[\log(1+\omega R_{t}Y)].
\]
Substituting $t=(1-{\rme}^{-s})/\theta=s^{*}$, $p=\theta s^{*}$,
yields a $\operatorname{GGC}(p,R_{s*}Y)$ variable, which is also a $\operatorname{GGC}(1,R_{s*}Y\xi
_{p})$. Statement (i) is straightforward. Statement (ii) is also
not difficult to verify.
\end{pf}

In order to evaluate the price in (\ref{Timeprice}), it remains to
evaluate $M_{1,\tau^{*}}$. We sketch out the details to do this
via the Double CFTP perfect sampler in \cite{DevroyeJames}. The
deterministic time change allows us to exploit generally the most
efficient case, $\theta=1$, of the Double CFTP.

First, note again that,
$R_{\tau^{*}}\stackrel{d}{=}k(\tau^{*},\tilde{U}\tau^{*})$ where
$k(t,y)$ is a known function, and $\tilde{U}$ is a Uniform$[0,1]$
random variable. Hence, in order to sample $R_{\tau^{*}}$ we simply
need to draw $\tilde{U}$. Note that we write $\tilde{U}$ to
distinguish it from the uniform variables we introduce below
denoted as $U_{i}$. Assuming
\[
D\stackrel{d}{=}R_{\tau^{*}}Y\xi_{p}\stackrel{d}{=}k(\tau
^{*},\tilde{U}\tau^{*})Y\xi_{p}
\]
is bounded by a positive constant $c$, the Double CFTP exact
sampler in \cite{DevroyeJames} is based on the following steps:
\begin{longlist}
\item[\textit{Backward phase}.] For $i=-1,-2,\ldots$: keep generating
$(U_i,D_i,D'_i)$ and
storing $(D_i,D'_i)$ until $U_\mathbb{T}\leq
|D_\mathbb{T}-D'_\mathbb{T}|/(2c)$. Keep $\mathbb{T}$.
\item[\textit{Set starting point}.] Set $M_{1,\tau^{*}}=D_\mathbb{T}\wedge
D'_\mathbb{T}+2cU_\mathbb{T}$.
\item[\textit{Forward phase}.] For $i=\mathbb{T}+1,\mathbb{T}+2,\ldots,-1$:
given $(D_i,D'_i,M_{1,\tau^{*}})$ previously stored, do the following
step: generate $U'$ uniform $[0,1]$, $\xi_{1/2}$, and generate $U$
uniform $[0,1]$, and construct $X=(1-U)M_{1,\tau^{*}}+UD_i\xi
_{1/2}+UD'_i(1-\xi_{1/2})$. Repeat this step until:
\[
U'\biggl[\mathbb{I}_{[0,1]}\biggl(\frac{X-M_{1,\tau^{*}}}{D_i-M_{1,\tau
^{*}}}\biggr)\frac{1}{|D_i-M_{1,\tau^{*}}|}
+\mathbb{I}_{[0,1]}\biggl(\frac{X-M_{1,\tau^{*}}}{D'_i-M_{1,\tau
^{*}}}\biggr)\frac{1}{|D'_i-M_{1,\tau^{*}}|}\biggr]> 1/c
\]
or $X< D_i\wedge D'_i$ or $X> D_i\vee D'_i$. Then set $M_{1,\tau^{*}}=X$.
\item[\textit{Output}.] Return $M_{1,\tau^{*}}$.
\end{longlist}
See \cite{DevroyeJames} for more details.
\begin{rem}
These results, which are considerably simplified by using the
deterministic time change, apply to a wide choice of kernels.
It would also be
nice to find models for $T(\tau)$ whose marginal distributions
were not strongly dependent on the form of the kernel. Even
better, would be the ability to specify laws in a manner similar
to how one selects the BDLP of an OU to induce general
self-decomposable laws for the instantaneous volatility $v(t)$. In
the next few sections, we will show that this can be done for
convoluted subordinators we refer to as quantile
clocks.
\end{rem}

%s3 ###
\section{Quantile clocks}\label{sec3}

As in the \hyperref[intro]{Introduction},
let $R$ denote a positive random variable with continuous strictly
increasing cumulative distribution function $F_{R}$. Let $Q_{R}$
denote its corresponding quantile function, that is, the continuous
inverse of the cumulative distribution function. Furthermore,
assume that $\mathbb{E}[R]<\infty$. Then for a subordinator $L$,
we say that the process $T_{R}=(T_{R}(t)\dvtx t\ge0)$, defined
as
%
%e3.1 ###
\begin{equation}
\label{defclock}
T_{R}(t):=\int_{0}^{t}Q_{R}\biggl(\biggl(1-\frac{s}{t}\biggr)_{+}\biggr)L(ds)\qquad
\mbox{for } t\ge0
\end{equation}
is a \textit{quantile clock} with
parameters $(R,L)$. Note here that $R$ does not depend on $t$.
Furthermore, $Q_{R}$ can be evaluated numerically in many cases, even
though it may not have a closed form.
\begin{prop}\label{Qclock1} A quantile clock $T_{R}=(T_{R}(t)\dvtx t\ge0)$
with parameters $(R,L)$ has the following properties:
\begin{longlist}[(iii)]
\item[(i)]If the support of the density of $R$, say $f_{R}$, is of
the form $[0,b), b>0$, that is, $Q_R(0)=0$, then
$T_{R}$ are random processes with samples paths that are almost sure
strictly continuous and strictly increasing.
\item[(ii)]Suppose the density of $R$ has support starting at $a>0$,
that is, $Q_R(0)=a$, then there is a positive random variable $\tilde
{R}$ with $Q_{\tilde{R}}(0)=0$, such that $R\stackrel{d}{=}\tilde
{R}+a$ and $Q_R(u)=Q_{\tilde{R}}(u)+a$ for $u\in[0,1]$. Hence, it
follows that the clock can be represented as
%
%e3.2 ###
\begin{equation}
\label{mixedclock} T_{R}(t)=T_{\tilde{R}}(t)+aL(t),\qquad t\ge0,
\end{equation}
where $T_{\tilde{R}}$ satisfies \textup{(i)}. Note $T_{\tilde{R}}$ is an
$(\tilde{R},L)$ quantile clock and is obviously not independent of $L$.
\item[(iii)]For each fixed $t$, the marginal distribution
\[
T_{R}(t)\stackrel{d}{=}\zeta(t),
\]
where $\zeta$ is a subordinator such that
$\zeta(1)$ is a random variable with Laplace exponent
\[
\psi_{\zeta}(\omega)=\mathbb{E}[\psi_{L}(\omega
R)]=\psi_{T_{R}(1)}(\omega).
\]
\item[(iv)]That is, the L\'{e}vy density of $\zeta$ has the form
\[
\rho_{\zeta}(s)=\int_{0}^{\infty}\rho_{L}(s/r)r^{-1}F_{R}(dr).
\]
Note furthermore that for a constant $c$,
$T_{cR}(t)\stackrel{d}{=}c\zeta(t)$.
\end{longlist}
\end{prop}
\begin{pf} Setting $k(t,s)=Q_{R}((1-\frac{s}{t})_{+}))$ it
follows from \cite{Bender}, Proposition 1, that in order to verify
statement (i) we only need to
check whether $k(t,s)$ satisfies conditions (a), (b), (c). Conditions (b)
and (c) are obvious and it remains to check the integrability
condition, which follows from
\[
\int_{0}^{t}Q_{R}\biggl(\biggl(1-\frac{s}{t}\biggr)_{+}\biggr)\,ds=t\int
_{0}^{1}Q_{R}(u)\,du=t\mathbb{E}[R]<\infty,
\]
since $Q_{R}(U)\stackrel{d}{=}R$. Now using this, and standard
results for linear functionals of L\'{e}vy processes, we see that
for each fixed $t$, the Laplace exponent of $T_{R}(t)$ is given by
\[
\int_{0}^{t}\psi_{L}\biggl(\omega
Q_{R}\biggl(\biggl(1-\frac{s}{t}\biggr)_{+}\biggr)\biggr)\,ds=t\int_{0}^{1}\psi_{L}(\omega
Q_{R}(u))\,du=t\mathbb{E}[\psi_{L}(\omega R)]
\]
verifying (iii). Statements (ii) and (iv) follows easily from (i) and
(iii). Note that the quantile function of $R$ in statement (ii)
violates condition (b).
\end{pf}

We now give an interesting example that has explicit laws.
\begin{example}[(Arcsine/Bessel occupation time quantile clocks driven
by $L=\gamma_{\theta}$)]\label{Arcsineexample} First, recall the
exponentially tilted stable subordinator $\widehat{S}_{\alpha}$ discussed
in Remark \ref{stableremark}.
Suppose that one specifies $R\stackrel{d}{=}\beta_{1/2,1/2}$ and
$L=\gamma_{\theta}$, a $\operatorname{gamma}(\theta)$ subordinator.
Then
\[
Q_{(\beta_{1/2,1/2})}(u)=\sin^{2}\biggl(\frac{\pi}{2}u\biggr),\qquad 0<u<1,
\]
and the quantile clock is defined as
\[
T_{\beta_{1/2,1/2}}(t)=\int_{0}^{t}
\sin^{2}\biggl(\frac{\pi(t-s)_{+}}{2t}\biggr)\gamma_{\theta}(ds)
\]
for $t\ge0$. It follows (see \cite{CifarelliMelilli,JamesLamperti,JRY}) that for
each fixed $t$,
\[
T_{\beta_{1/2,1/2}}(t)\stackrel{d}{=}\widehat{S}_{1/2}(\gamma_{2\theta
t}/2)\stackrel{d}{=}\gamma_{\theta t}\beta_{\theta t+1/2,\theta
t+1/2}.
\]
More generally, for each fixed $0<\alpha<1$, let
$\mathbb{O}_{\alpha}(s)=\int_{0}^{s}\indic_{(B_{u}>0)}\,du$ denote
the time spent
positive up till time $s$ of a symmetrized Bessel process
$(B_{u}, u\ge0)$ of dimension $2-2\alpha$, see \cite{BPY}. Then setting
$R\stackrel{d}{=}\mathbb{O}_{\alpha}(1):=\mathbb{O}_{\alpha}$, the
quantile of $\mathbb{O}_{\alpha}$ is
\[
Q_{\mathbb{O}_{\alpha}}(u)=\frac{Q_{X_{\alpha}}(u)}{Q_{X_{\alpha
}}(u)+1}\qquad\mbox{where }Q_{X_{\alpha}}(u)={\biggl[\frac{\sin(\pi\alpha
u)}{\sin(\pi\alpha(1-u))}\biggr]}^{1/\alpha}
\]
is the quantile function of the ratio of i.i.d. positive stable
random variables $X_{\alpha}=S_{\alpha}/S'_{\alpha}$. Then, from
James \cite{JamesLamperti} (see Section 7), the clock $T_{\mathbb
{O}_{\alpha}}$
with parameters $(\mathbb{O}_{\alpha},\gamma_{\theta})$, satisfies
for each
fixed $t$,
\[
T_{\mathbb{O}_{\alpha}}(t)\stackrel{d}{=}\widehat{S}_{\alpha}(\gamma
_{\theta
t/\alpha}/2)\stackrel{d}{=}\gamma_{\theta t}\mathbb{O}_{\alpha
,\theta t},
\]
where
\begin{eqnarray*}
\mathbb{O}_{\alpha,\theta t}&\stackrel{d}{=}&\beta_{\theta
t,1}\mathbb{O}_{\alpha,\theta t}+(1-\beta_{\theta t,1})\mathbb
{O}_{\alpha}\\
&\stackrel{d}{=}&\beta_{\alpha+\theta t,1-\alpha}\mathbb
{O}_{\alpha,\alpha+\theta t}+(1-\beta_{\alpha+\theta t,1-\alpha
})\xi_{1/2}
\end{eqnarray*}
are random variables corresponding to the time spent positive of
generalized Bessel bridges as explained in \cite{JamesLamperti},
Section 5. These variables can be exactly sampled in various
ways as explained in \cite{DevroyeJames}. Furthermore, from
\cite{JamesLamperti}, Proposition 5.3, it follows that for $0<p=\theta t\leq1$,
\[
T_{\mathbb{O}_{\alpha}}(p/\theta)\stackrel{d}{=}\gamma
_{1}\tilde{\mathbb{O}}_{\alpha,p},
\]
where $\tilde{\mathbb{O}}_{\alpha,p}\stackrel{d}{=}\beta
_{p,1-p}\mathbb{O}_{\alpha,p}$ is $\operatorname{GGC}(1, \mathbb{O}_{\alpha}\xi
_{p})$ with density
%
%e3.3 ###
\begin{eqnarray}\label{denB}
f_{\tilde{\mathbb{O}}_{\alpha,p}}(y)&=&\frac{2^{{p}/{\alpha
}}}{\pi}
{y^{p-1}\sin\biggl(
\frac{p}{\alpha}\arctan\biggl(\frac{{(1-y)}^{\alpha}\sin(\pi
\alpha)}{{(1-y)}^{\alpha}\cos(\pi\alpha)+y^{\alpha}}\biggr)
\biggr)}\nonumber\\[-8pt]\\[-8pt]
&&{}\times [y^{2\alpha}+2y^{\alpha}{(1-y)}^{\alpha}\cos(\alpha
\pi)+{(1-y)}^{2\alpha}]^{{-p}/({2\alpha})},\nonumber
\end{eqnarray}
$0<y<1$. In general, the process
$(\widehat{W}_{\mu}(T_{\mathbb{O}_{\alpha}}(t)),t \ge0)$ has
almost surely
continuous sample paths and satisfies, for each fixed $t$,
\[
\mathbb{E}\bigl[{\rme}^{i\omega\widehat{W}_{\mu}(T_{\mathbb{O}_{\alpha
}}(t))}\bigr]=2^{\theta
t/\alpha}\bigl(1+{\bigl(1+(\omega^{2}/2-i\mu\omega)\bigr)}^{\alpha}\bigr)^{-\theta
t/\alpha}.
\]
\end{example}
\begin{rem}The last example shows that the class of quantile
clocks where $L$ is a $\operatorname{GGC}(\theta,Y)$ subordinator is equivalent
in a marginal sense to the representation of GGC variables in
terms of Wiener--Gamma integrals as defined and presented in
\cite{JRY}. That manuscript, along with the works of
\cite{BFRY,JamesBernoulli,JamesLamperti,JLP}, yield many
examples of quantile functions which can be used to construct
quantile clocks with explicit laws, of which quite a few are
constructed from $Q_{X_{\alpha}}$. We shall encounter some more
examples in Section \ref{sec6}, although in a slightly different context.
\end{rem}

%s4 ###
\section{Choosing quantile clocks to have specific laws in
$\mathcal{U_{\delta}}$}\label{sec4}

The results in the previous section
suggest that the the marginal distributions of the quantile
clocks, while equating nicely to the marginals of a subordinator
$\zeta$, are strongly dependent on a random variable $R$, induced by $Q_{R}$.
Noting that $Q_{R}(u)$ is in fact a deterministic function one
would like to be able to choose explicit laws for
$T_{R}$, regardless of the function $Q_{R}$. That is to say, how
does one choose $L$ so that $T_{R}(t)$ has a marginal distribution
not obviously depending on $R?$ For example, for each $Q_{R}$ how
does one choose $L$ so that
$T_{R}(t)\stackrel{d}{=}\widehat{S}_{\alpha}(t)?$ Or how does one choose
$T_{R}(t)$ so that a log price process $\widehat{W}_{\mu}(T_{R}(t))$
has marginal distributions that are equivalent to a CGMY
process? Finally, for different quantile functions $Q_{R_{1}}$,
$Q_{R_{2}}$, $R_{1}$ not equivalent\vspace*{-1pt} in distribution to $R_{2}$,
how to choose the driving L\'{e}vy processes, say $L_{1}$ and
$L_{2}$, such that marginally for each fixed $t$,
$T_{R_{1}}(t)\stackrel{d}{=}T_{R_{2}}(t)?$

We saw that this was difficult in the case of general convoluted
subordinators as their laws depend strongly on $t$ through the
kernel or variable $k(t,Ut)$. However, Proposition
\ref{Qclock1} shows that one can represent
%
%e4.1 ###
\begin{equation}
\label{TRV} T_{R}(t)\stackrel{d}{=}\int_{0}^{1}Q_{R}(y)L(t\,dy),
\end{equation}
and there is a clear separation of the effects of $t$ and $Q_{R}$.
This is similar to the case of the OU models $v(0)$, see
\cite{BN98,BNS2001,JV}, where every positive self-decomposable
random variable can be represented as
%
%e4.2 ###
\begin{equation}
\label{SD1}
v(0)\stackrel{d}{=}v(t)\stackrel{d}{=}\int_{-\infty}^{y}{\rme}^{-\lambda(y-s)}\vartheta(\lambda
\,ds)\stackrel{d}{=}\int_{0}^{\infty}{\rme}^{-s}\vartheta(ds),
\end{equation}
where $\vartheta$, is a subordinator referred to as a OU--BDLP. More
strikingly, there is a simple way of obtaining any desired
self-decomposable law for $v(0)$ by choosing the BDLP according to
either of the equations
%
%e4.3 ###
\begin{equation}
\label{SD2}
\psi_{\vartheta}(\omega)=\omega\psi'_{v(0)}(\omega)
\quad\mbox{and}\quad\rho_{\vartheta}(x)=-\rho_{v(0)}(x)-x\rho'_{v(0)}(x).
\end{equation}

We noticed from (\ref{TRV}) that if
$R\stackrel{d}{=}U^{1/\delta}\stackrel{d}{=}\beta_{\delta,1}$ for
$\delta>0$, then
\[
T_{\beta_{\delta,1}}(t)\stackrel{d}{=}\int_{0}^{1}u^{1/\delta
}Z(t\,du)\stackrel{d}{=}\zeta^{(\delta)}(t),
\]
where, we substitute $Z$ for $L$, and $\zeta^{(\delta)}$ are subordinators
having laws in Jurek's \cite{IksanovJurek,Jurek85,Jurek88,Jurek89}
$\mathcal{U}_{\delta}$ class of
\textit{generalized} $s$-selfdecomposbale laws, where
$\mathcal{U}_{\delta}\subset\mathcal{J}$. The case of $\delta=1$,
corresponds to Jurek's $\mathcal{U}=\mathcal{U}_{1}$ class of
$s$-selfdecomposable class. Using Jurek \cite{Jurek88,Jurek89}, one
sees that for $0<\delta_{1}<1<\delta_{2}<\infty$,
\[
\mathcal{G}_{+}\subset\mathcal{G} \subset\mathcal{L} \subset
\mathcal{U}_{\delta_{1}} \subset\mathcal{U}\subset
\mathcal{U}_{\delta_{2}} \subset\mathcal{J}.
\]
It follows that for each $\zeta^{(\delta)} \in\mathcal{U}_{\delta}$
there is a, $\mathcal{U}_{\delta}$-BDLP, $Z$ such that
\[
\psi_{\zeta^{(\delta)}}(\omega)=\psi_{T_{\beta_{\delta
,1}}(1)}(\omega)=\int_{0}^{1}\psi_{Z}(\omega
u^{1/\delta})\,du=\omega^{-\delta}\int_{0}^{\omega}\psi_{Z}(
u)\delta u^{\delta-1}\,du
\]
and hence from \cite{CJurek}, Lemma 1, which can be verified
directly by taking derivatives with respect to $\omega$ of
%
%e4.4 ###
\begin{equation}
\label{firstderivative}
\omega^{\delta}\psi_{\zeta^{(\delta)}}(\omega)=\int_{0}^{\omega
}\psi_{Z}(
u)\delta u^{\delta-1}\,du,
\end{equation}
one sees that the $\mathcal{U}_{\delta}$-BDLP $Z$ is related to
$\zeta^{(\delta)}$, and hence $T_{\beta_{\delta,1}}$, by the
equation
%
%e4.5 ###
\begin{equation}
\label{Uvanilla}
\psi_{Z}(\omega)=\psi_{\zeta^{(\delta)}}(\omega)+\frac{1}{\delta}
\omega\psi'_{\zeta^{(\delta)}}(\omega).
\end{equation}
This is
analogous to the relationships between $v(0)$ and its OU--BDLP
$\vartheta$, given in (\ref{SD1}) and (\ref{SD2}). We shall show
that this relationship becomes more explicit as one restricts
their choices of laws for $\zeta^{(\delta)}$ to $\mathcal{L}$,
$\mathcal{G}$ and $\mathcal{G}_{+}$.
\begin{rem}
The specifications in (\ref{Uvanilla}) and its refinements now
allow us to specify any law in $\mathcal{U}_{\delta}$ for quantile
clocks based on $Q_{U^{1/\delta}}$, analogous to the case of the BNS--OU
$v(t)$. This, as far as we know, is the first instance where such
a property has been noticed for convoluted subordinators. However,
in terms of choices of $Q_{R}$ this is still restrictive.
The next results show how, for a large class of quantile
functions $Q_{R}$, to choose $L$ such that for each fixed $t$,
$T_{R}(t)\stackrel{d}{=}\zeta^{(\delta)}(t) \in\mathcal{U}_{\delta}$.
\end{rem}
\begin{theorem}\label{TheoremRY}Consider the specifications for a
quantile clock
\[
T_{R}(t)=\int_{0}^{t}Q_{R}\biggl(\biggl(1-\frac{s}{t}\biggr)_{+}\biggr)L(ds)
\]
with parameters
$(R,L)$. Now select $R$ so that its density has bounded support
and let $Y$ denote a positive bounded random variable such that
%
%e4.6 ###
\begin{equation}\label{betacondition}
RY\stackrel{d}{=}U^{1/\delta}\stackrel{d}{=}\beta_{\delta,1}
\end{equation}
for a fixed $\delta>0;$ if $R\stackrel{d}{=}cU^{1/\delta}$, then
$Y=1/c$. Suppose that one wishes to choose $L$ such that
\[
T_{R}(t)\stackrel{d}{=}\zeta^{(\delta)}(t),
\]
where
$\zeta^{(\delta)}$ is a subordinator with
\[
\psi_{\zeta^{(\delta)}}(\omega)=\mathbb{E}[\psi_{Z}(\omega
U^{1/\delta})],
\]
where $Z$ is a $\mathcal{U}_{\delta}$-BDLP, satisfying
(\ref{Uvanilla}), and hence $\zeta^{(\delta)}$ is in $\mathcal
{U}_{\delta}$. Then for this $Z$, $L$~is chosen such
that
%
%e4.7 ###
\begin{equation}
\label{Lspec} \psi_{L}(\omega)=\mathbb{E}[\psi_{Z}(\omega Y)]\quad
\mbox{equivalently}\quad
\rho_{L}(x)=\mathbb{E}[\rho_{Z}(x/Y)Y^{-1}].
\end{equation}
That is,
%
%e4.8 ###
\begin{equation}
\label{genLU}
\psi_{L}(\omega)=\mathbb{E}\bigl[\psi_{\zeta^{(\delta)}}(\omega
Y)\bigr]+\frac{1}{\delta} \omega
\mathbb{E}\bigl[Y\psi'_{\zeta^{(\delta)}}(\omega Y)\bigr].
\end{equation}
Note that $Y$ is chosen independent of $R$ and $Z$.
\end{theorem}
\begin{pf}The difficulty of this result is envisioning its construction.
The proof itself is otherwise straightforward, since (\ref{Lspec}) and
(\ref{Uvanilla}) implies that
\[
\psi_{\zeta^{(\delta)}}(\omega)=\mathbb{E}[\psi_{L}(\omega
R)]=\mathbb{E}[\psi_{Z}(\omega U^{1/\delta})].
\]
\upqed\end{pf}

We now specialize this result to self-decomposable laws.
\begin{theorem}\label{SelfDecomposableThm} Consider quantile clocks
$T_{R}$ with parameters $(R,L)$ satisfying (\ref{betacondition}) and
(\ref{Lspec}), (\ref{genLU}). The next result describes further
specifications in order for
$T_{R},(\zeta^{(\delta)})$ to have laws in
$\mathcal{L},\mathcal{G}$ and $\mathcal{G}_{+}$, respectively.
\begin{enumerate}[III.]
\item[I.] $T_{R}\in\mathcal{L}$: If $T_{R}$ is selected such
that its marginal laws are self-decomposable, then it is known
that there exists a subordinator $\vartheta$, such that
\[
T_{R}(1)\stackrel{d}{=}\zeta^{(\delta)}(1)\stackrel{d}{=}\int
_{-\infty}^{y}{\rme}^{-\lambda(y-s)}\vartheta(\lambda
\,ds)\stackrel{d}{=}\int_{0}^{\infty}{\rme}^{-s}\vartheta(ds)\stackrel{d}{=}v(0).
\]
Furthermore, adapting (\ref{SD2}), one has
%
%e4.9 ###
\begin{equation}
\label{OUcondition2}
\psi_{\vartheta}(\omega)=\omega\psi'_{\zeta^{(\delta)}}(\omega
)\quad\mbox{and}\quad\rho_{\vartheta}(x)=-\rho_{\zeta^{(\delta)}}(x)-x\rho'_{\zeta
^{(\delta)}}(x).
\end{equation}
\begin{enumerate}[(ii)]
\item[(i)]Hence, the BDLP $L$ has to be chosen such
that the L\'{e}vy density of $Z$ is
\begin{eqnarray*}
\rho_{Z}(x) &=& \biggl(1-\frac{1}{\delta}\biggr)\rho_{\zeta^{(\delta
)}}(x)-\frac{x}{\delta}\rho'_{\zeta^{(\delta)}}(x) \\
&=&
\rho_{\zeta^{(\delta)}}(x)+\frac{1}{\delta}\rho_{\vartheta}(x).
\end{eqnarray*}
That is,
\[
\rho_{L}(x)=\mathbb{E}
\biggl[\biggl(\rho_{\zeta^{(\delta)}}(x/Y)+\frac{1}{\delta}\rho_{\vartheta
}(x/Y)\biggr)Y^{-1}\biggr].
\]
\item[(ii)]Statement \textup{(i)} implies that the subordinators are related as
follows:
\[
Z(t)
\stackrel{d}{=}\zeta^{(\delta)}(t)+\vartheta(t/\delta),\qquad t\ge0.
\]
\end{enumerate}
\item[II.] $T_{R} \in\mathcal{G}$: If $T_{R},(\zeta^{(\delta)})$ is
selected such
that its marginal laws are $\operatorname{GGC}(\nu)$, it follows that
\[
\rho_{\zeta^{(\delta)}}(x)=x^{-1}\int_{0}^{\infty}{\rme}^{-xy}\nu(dy)
\quad\mbox{and}\quad\rho_{\vartheta}(x)
=\int_{0}^{\infty}{\rme}^{-xy}y\nu(dy).
\]
\begin{enumerate}[(ii)]
\item[(i)]Hence, the L\'{e}vy density of $Z$, say $\rho_{Z}$, satisfies
\[
\rho_{Z}(x)=\int_{0}^{\infty}{\rme}^{-xy}[x^{-1}+y/\delta
]\nu(dy).
\]
\item[(ii)]Equivalently,
\[
\psi_{Z}(\omega)=\int_{0}^{\infty}\biggl[\log(1+\omega/y)+\frac
{1}{\delta}\frac{\omega}{y+\omega}\biggr]\nu(dy).
\]
\end{enumerate}
\item[III.] $T_{R} \in\mathcal{G}_{+}$: If $T_{R},(\zeta^{(\delta
)})$ is selected such
that its marginal laws are $\operatorname{GGC}(\theta,V)$,

\begin{enumerate}[(ii)]
\item[(i)]then $L$ must be selected such that it is equivalent in
distribution to the subordinator
\[
L(s)\stackrel{d}{=}\zeta_{\delta,Y}(s)+\vartheta_{Y}(s\theta
/\delta),\qquad s\ge0,
\]
where $\zeta_{\delta,Y}$ is a $\operatorname{GGC}(\theta, VY)$ subordinator and
\[
\vartheta_{Y}(s)\stackrel{d}{=}\sum_{k=1}^{N(s
)}\gamma^{(k)}_{1}V_{k}Y_{k},\qquad
s\ge0,
\]
where $(\gamma^{(k)}_{1})$ are independent $\operatorname{exponential}(1)$ variables,
$(Y_{k})$ are i.i.d. variables with distribution $F_{Y}$, $(V_{k})$
are i.i.d. $F_{V}$, and $N(s)$ denotes a homogeneous Poisson process
with $\mathbb{E}[N(s)]=s$.
\item[(ii)]As special cases
$T_{R}(t)\stackrel{d}{=}\gamma_{\theta}(t)$ is obtained by setting
$V=1$.
\end{enumerate}
\end{enumerate}
\end{theorem}

The proof of this result is fairly immediate from the definitions of
the various classes, details are omitted.
\begin{rem}Theorem \ref{SelfDecomposableThm} shows that in order to
specify $T_{R}$ to have laws in $\mathcal{L}$, one only needs to
identify the OU--BDLP $\vartheta$ that leads to a corresponding
stationary law for $v(t)\stackrel{d}{=}T_{R}(1)\stackrel{d}{=}\zeta
^{(\delta)}(1)$, and use it appropriately to define $L$. One may
consult for instance \cite{BNS2001} for many explicit examples
$\vartheta$, and the laws they induce.
\end{rem}

%s5 ###
\section{Choosing $R$ and $Y$ such that $RY\stackrel{d}{=}U^{1/\delta}$}\label{sec5}

The results in the previous section show that for a deterministic
quantile function $Q_{R}$ one can choose quite arbitrary marginal
laws for $T_{R}$, analogous to the case of $v(0)$, provided that
one identifies a variable $Y$ such that
$RY\stackrel{d}{=}U^{1/\delta}=\beta_{\delta,1}$ for a fixed
$\delta>0$. Notice that
%
%e5.1 ###
\begin{equation}
\label{genholmgren} Q_{R^{1/\delta}}(u)={[Q_{R}(u)]}^{1/\delta}.
\end{equation}

The easiest case is to choose $R=U$ and $Y=1$, which as seen from
 (\ref{genholmgren}) leads to quantile clocks corresponding to the
\textit{Holmgren--Liouville} convoluted subordinators discussed in
\cite{Bender}. The equation (\ref{genholmgren}) suggests that one
may always work with the pair satisfying the solution $RY=U$ and
then obviously $R^{1/\delta}Y^{1/\delta}=U^{1/\delta}$. However,
the case of $\delta=1$ may not always be the most obvious.
\begin{example}[(Beta variables including the arcsine
distribution)]\label{betaexample}
Consider the case of products of independent beta variables
%
%e5.2 ###
\begin{equation}
\label{BArcsine}
\beta_{\delta,\kappa-\delta}\beta_{\kappa,1+\delta-\kappa
}\stackrel{d}{=}\beta_{\delta,1}\stackrel{d}{=}U^{1/\delta}
\end{equation}
for $\delta\leq\kappa\leq1+\delta$. Hence, for each $\delta$, one
can choose many $(R,Y)$, ranging over $\delta\leq\kappa\leq
1+\delta$, such that
%
%e5.3 ###
\begin{eqnarray}
\label{betachoice}
(R^{1/\delta},Y^{1/\delta})&\stackrel{d}{=}&(\beta_{\delta,\kappa
-\delta},\beta_{\kappa,1+\delta-\kappa})
\quad\mbox{or}\nonumber\\[-8pt]\\[-8pt]
(R^{1/\delta},Y^{1/\delta})&=&(\beta_{\kappa,1+\delta-\kappa
},\beta_{\delta,\kappa-\delta}).\nonumber
\end{eqnarray}
Furthermore, for some $b>0$, and each fixed $\delta$ the variables
in (\ref{betachoice}) lead to variables $R^{1/b}$ and $Y^{1/b}$,
not having beta distributions, that satisfy
$(RY)^{1/b}\stackrel{d}{=}U^{1/b}$.
\end{example}

Lets look at a special case of this in more detail.
\begin{example}[(Kumaraswamy and generalized arcsine clocks)]
Setting $\kappa=1$ and $\delta=\alpha$, (\ref{betachoice}) leads
to the choice of the pair
%
%e5.4 ###
\begin{equation}
\label{Arcpair}
(\beta_{\alpha,1-\alpha},1-U^{1/\alpha})\stackrel{d}{=}(\beta
_{\alpha,1-\alpha},\beta_{1,\alpha}),
\end{equation}
such that
$R^{1/\alpha}Y^{1/\alpha}\stackrel{d}{=}U^{1/\alpha}\stackrel
{d}{=}\beta_{\alpha,1}$,
where the first component in (\ref{Arcpair}) has the
\textit{generalized arcsine law} which arises in many studies of
random processes. Setting
$R^{1/b}\stackrel{d}{=}{[1-U^{1/\alpha}]}^{1/b}=K_{\alpha,b}$ leads
to the quantile function
\[
Q_{K_{\alpha,b}}(u)={[1-(1-u)^{1/\alpha}]}{}^{1/b}
\]
of a Kumaraswamy distribution. Hence, the law of a Kumaraswamy
quantile clock $T_{K_{\alpha,b}}$, that is, with parameters
$(K_{\alpha,b},L)$, can be specified such that its marginals
satisfy
\[
T_{K_{\alpha,b}}(t)\stackrel{d}{=}\zeta^{(\alpha b)}(t),
\]
where $\zeta^{(\alpha b)}$ is a subordinator having any law in
$\mathcal{U}_{\alpha b}$ and hence in $\mathcal{L}$. Specifically,
this is done by the choice of
\[
\psi_{L}(\omega)=\mathbb{E}\bigl[\psi_{\zeta^{(\alpha b)}}(\omega
{(\beta_{\alpha,1-\alpha})}^{1/b})\bigr]+\frac{1}{\alpha b} \omega
\mathbb{E}\bigl[\psi'_{\zeta^{(\alpha b)}}(\omega
{(\beta_{\alpha,1-\alpha})}^{1/b}){(\beta_{\alpha,1-\alpha})}^{1/b}\bigr].
\]
Note that if we instead choose $R=\beta_{1/2,1/2}$ and hence
$R^{1/b}\stackrel{d}{=}(\beta_{1/2,1/2})^{1/b}$ we obtain clocks
based on the arcsine law with quantiles
%
%e5.5 ###
\begin{equation}
\label{TrueArcsine}
Q_{R^{1/b}}(u):=\bigl[Q_{(\beta_{1/2,1/2})}(u)\bigr]^{1/b}=\sin^{2/b}
\biggl(\frac{\pi}{2}u\biggr).
\end{equation}
This case can be compared with Example \ref{Arcsineexample}.
\end{example}

%s5.1 ###
\subsection{Selections based on decompositions of an $\operatorname{exponential}(1)$ variable}\label{sec5.1}

In general, by rescaling to $[0,1]$, we see that
choosing an $R$ and $Y$ to satisfy (\ref{betacondition}) for some
$\delta>0$ is equivalent to choosing variables $\ell_{R}$ and
$\ell_{Y}$ such that
%
%e5.6 ###
\begin{equation}
\label{expcond} \ell_{R}+\ell_{Y}\stackrel{d}{=}\gamma_{1}/\delta.
\end{equation}

Recall also the relationship between quantiles of a positive
variable $X$ and ${\rme}^{-X}$,
\[
Q_{\rme^{-X}}(u)={\rme}^{-Q_{X}(1-u)},\qquad 0\leq u\leq1.
\]

There are obviously many pairs satisfying (\ref{expcond}). We next
look at two different types of examples based on suggestions made
to us by Prof. Marc Yor.
\begin{example}[{[Fractional and integer parts of an $\operatorname{exponential}(1)$]}]\label{exam5.3}
We first note that one of the reasons the following example is
interesting is that it identifies a concrete example of a quantile
clock that is of the form
\[
T_{R}(t)=T_{\tilde{R}}+aL(t)
\]
as specified in statement (ii) of Proposition \ref{Qclock1}, but
where we can apply Theorem~\ref{TheoremRY}. Now, following Chaumont and
Yor (\cite{Chaumont}, page 42, Exercise 2.18), let
$[\gamma_{1}]$ and
${\{\gamma_{1}\}}$ denote the fractional part and integer part of
an $\operatorname{exponential}(1)$ variable~$\gamma_{1}$, then (remarkably) these
variables are independent and obviously satisfy
\[
[\gamma_{1}]+{\{\gamma_{1}\}}=\gamma_{1}.
\]
In this case, $RY\stackrel{d}{=}U$, for
%
%e5.7 ###
\begin{equation}
\label{geometric}
R\stackrel{d}{=}{\rme}^{-[\gamma_{1}]}\stackrel{d}{=}U(1-e^{-1})+e^{-1}\quad
\mbox{and}\quad
Y\stackrel{d}{=}{\rme}^{-\{\gamma_{1}\}},
\end{equation}
where
$\{\gamma_{1}\}$ is a geometric random variable with success
probability $1-e^{-1}$ and values in $\{0,1,2,\ldots\}$. We say
such a variable is geometric $(1-e^{-1})$. We can extend this case as follows.
\begin{prop}For $0<p\leq1$, let $\tilde{U}_{p}\stackrel
{d}{=}Up+(1-p)$ and let $X_{p}$ be a $\operatorname{geometric}(p)$ variable. Then
\[
\tilde{U}_{p}{\rme}^{-X_{p}[-{\log}(1-p)]}\stackrel{d}{=}U.
\]
\end{prop}
\begin{pf}It is easy to verify that the Laplace transforms of $-{\log}
(Up+(1-p))$ and $X_{p}[-{\log}(1-p)]$ are given, respectively, by
\[
\frac{1-(1-p)^{(1+\omega)}}{p(1+\omega)}
\quad\mbox{and}\quad\frac
{p}{1-(1-p)^{(1+\omega)}}.
\]
Hence, their product is $1/(1+\omega)$, which is the desired result.
\end{pf}

That is, there are variables
$\tilde{U}_{p}Y_{p}\stackrel{d}{=}U$ for
\[
R=\tilde{U}_{p}\stackrel{d}{=}Up+(1-p)
\quad\mbox{and}\quad Y_{p}={\rme}^{-X_{p}[-{\log}(1-p)]}
\]
for $X_{p}$ a $\operatorname{geometric}(p)$ variable for each $0<p\leq1$.
Naturally $\tilde{U}_{1}\stackrel{d}{=}U$. Hence, for each fixed $p$,
\[
T_{\tilde{U}_{p}}(t)=L(t)(1-p)+p\int_{0}^{t}\biggl(1-\frac{s}{t}\biggr)_{+}L(ds)\qquad
\mbox{for } t\ge0.
\]
We now state a result which now follows obviously from Theorem \ref
{TheoremRY} and applies for quantiles clocks based on
the variable $\tilde{U}^{1/\delta}_{p}$.
\begin{prop}For each $\delta>0$, and $0<p\leq1$, set $\lambda=-{\log}
(1-p)$, then the quantile clock
\[
T_{\tilde{U}^{1/\delta}_{p}}(t)=\int_{0}^{t}\biggl[(1-p)+p\biggl(1-\frac
{s}{t}\biggr)_{+}\biggr]^{1/\delta}L(ds)
\]
can be specified such that for each $t$,
$T_{\tilde{U}^{1/\delta}_{p}}(t)\stackrel{d}{=}\zeta^{(\delta
)}(t)\in
\mathcal{U}_{\delta}$ if $L$ is chosen such that
\[
\psi_{L}(\omega)=\mathbb{E}\bigl[\psi_{\zeta^{(\delta)}}(\omega{\rme}^{-\lambda X_{p}/\delta})\bigr]+\frac{1}{\delta} \omega\mathbb
{E}\bigl[{\rme}^{-\lambda X_{p}/\delta}\psi'_{\zeta^{(\delta)}}(\omega{\rme}^{-\lambda X_{p}/\delta})\bigr],
\]
where $X_{p}$ is a $\operatorname{geometric}(p)$ random variable. When $p=1$, the
quantile clock is continuous, otherwise it has jumps.
\end{prop}
\end{example}
\begin{example}[(Splitting the Laplace exponent of $\gamma_{1}$, part
I)]\label{split}
Next, consider the Laplace exponent of $\gamma_{1}/\delta$,
\[
\log(1+\omega/\delta)=\int_{0}^{\infty}(1-{\rme}^{-s\omega/\delta})s^{-1}{\rme}^{-s}\,ds,
\]
then choose $\ell_{R}$ and $\ell_{Y}$ according to the
decomposition of the L\'{e}vy density
\[
s^{-1}{\rme}^{-s}=\pi_{1}(s)+\pi_{2}(s).
\]
In particular, $\ell_{R}$ and $\ell_{Y}$ are infinitely divisible
based on the L\'{e}vy densities $\pi_{1}$ and~$\pi_{2}$,
respectively. The simplest case is where
\[
s^{-1}{\rme}^{-s}=(1-\alpha)s^{-1}{\rme}^{-s}+\alpha
s^{-1}{\rme}^{-s}
\]
leading to
\[
(RY)^{1/\delta}\stackrel{d}{=}{[{\rme}^{-\gamma_{\alpha}}{\rme}^{-\gamma_{1-\alpha}}]}^{1/\delta}\stackrel{d}{=}U^{1/\delta}.
\]
So, for instance, the quantile clock with parameters $({\rme}^{-\gamma_{\alpha}},L)$, that is,
\[
T_{{\rme}^{-\gamma_{\alpha}}}(t)=\int_{0}^{t}{\rme}^{-Q_{\gamma_{\alpha}}(1-(1-s/t)_{+})}L(ds)
\]
is based on a nontrivial quantile which can however be evaluated
by various computational packages. Furthermore, our results show
that despite the complexity of this clock we can choose quite
general marginal laws for $T_{{\rme}^{-\gamma_{\alpha}}}$, not
directly depending on the quantile function, by working with a
BDLP satisfying
\[
\psi_{L}(\omega)=\mathbb{E}[\psi_{Z}(\omega{\rme}^{-\gamma_{1-\alpha}})]=\int_{0}^{1}\psi_{Z}(\omega
y)\frac{{[-{\log}(y)]}^{-\alpha}}{\Gamma(1-\alpha)}\,dy.
\]
\end{example}

%s6 ###
\section{GGC decompositions of a $\gamma$ subordinator and resulting clocks}\label{sec6}

We now identify a very large class of variables satisfying
(\ref{expcond}) based on Example \ref{split}, using variables in
$\mathcal{B}$ and in particular, the GGC class $\mathcal{G}_{+}$.
\begin{theorem}\label{decomposition}Let $\Omega_{\delta}$ be a
subordinator in $\mathcal{B}$ with L\'{e}vy density $\rho_{\Omega_{
\delta}}(s)=\int_{\delta}^{\infty}{\rme}^{-sx}q(x)\,dx$, where
$q(x)$ is a nonnegative measure such that $q(x)\leq1$ for \mbox{$x\ge
\delta$}. Then, there is another subordinator in $\mathcal{B}$, say
$\widehat{\Omega}_{\delta}$, with L\'{e}vy density $\rho_{\widehat
{\Omega}_{
\delta}}(s)=\int_{\delta}^{\infty}{\rme}^{-sx}[1-q(x)]\,dx$, such
that the L\'{e}vy density of a $\operatorname{gamma}(1)$
subordinator with scale $1/\delta$, has the decomposition
%
%e6.1 ###
\begin{equation}
\label{gammalevydecomp}
s^{-1}{\rme}^{-s\delta}=\rho_{\widehat{\Omega}_{ \delta}}(s)+\rho_{\Omega_{
\delta}}(s).
\end{equation}
Hence, the gamma subordinator can be expressed as a sum of the
subordinators $\Omega_{\delta}$ and $\widehat{\Omega}_{\delta}$,
which implies for each fixed $\theta$
\[
\gamma_{\theta}/\delta\stackrel{d}{=}\widehat{\Omega}_{\delta
}(\theta)+\Omega_{\delta}(\theta).
\]
\end{theorem}
\begin{pf}The L\'{e}vy density of the subordinator
$(\gamma_{1}(t)/\delta\dvtx t\ge0)$ can be expressed as
\[
s^{-1}{\rme}^{-s\delta}=\int_{\delta}^{\infty}{\rme}^{-sx}\,dx,
\]
leading easily to (\ref{gammalevydecomp}).
\end{pf}

We next describe an interesting special case involving variables
in $\mathcal{G}_{+}$.
\begin{theorem}\label{decompositionG}
Assume that $\Omega_{\delta}$ in Theorem \ref{decomposition} is a
$\operatorname{GGC}(1,V/\delta)$ subordinator for $V$ a random variable in
$[0,1]$, and let $X\stackrel{d}{=}(1-V)/V$. Let $(\Sigma_{t}(V)\dvtx t\ge
0)$ denote a subordinator with L\'{e}vy density denoted as
$\rho_{\Sigma_{1}(V)}$. Note that $\Sigma_{1}(V)$ is not random in
$V$.
\begin{longlist}[(iii)]
\item[(i)]The L\'{e}vy density of a $\operatorname{gamma}(1)$ subordinator with scale
$1/\delta$ has the decomposition
\[
s^{-1}{\rme}^{-s\delta}=\delta\rho_{\Sigma_{1}(V)}(s\delta)
+s^{-1}\mathbb{E}[{\rme}^{-s\delta/V}],
\]
where
%
%e6.2 ###
\begin{eqnarray}
\label{straddle}
\rho_{\Sigma_{1}(V)}(s)&=&\int_{1}^{\infty}{\rme}^{-sx}[1-F_{1/V}(x)]\,dx\nonumber\\[-8pt]\\[-8pt]
&=&\frac{1}{s}{\rme}^{-s}(1-\mathbb{E}[{\rme}^{-sX}]).\nonumber
\end{eqnarray}
\item[(ii)]If $0<\varrho=\mathbb{E}[-{\log}(V)]<\infty$, then
$\rho_{\Sigma_{1}(V)}(s)=\varrho f_{\Delta_{\mathbf{e}}(X)}(s)$,
where\break
$f_{\Delta_{\mathbf{e}}(X)}(s)$ is a density of a random variable
denoted as $\Delta_{\mathbf{e}}(X)$, determined
by~(\ref{straddle}). In this case $\Sigma_{t}(V)$ is a compound
Poisson process representable as
\[
\Sigma_{t}(V)=\sum_{k=1}^{N(\varrho t)}\Delta_{k},\qquad t\ge0.
\]
$(\Delta_{k})$ are i.i.d. random variables equal in distribution to
$\Delta_{\mathbf{e}}(X)$.
\item[(iii)]In general, for each fixed $\theta$,
%
%e6.3 ###
\begin{equation}
\label{GGCdecomp}
\gamma_{\theta}\stackrel{d}{=}\Sigma_{
\theta}(V)+\gamma_{\theta}M_{\theta},
\end{equation}
where
$M_{\theta}\stackrel{d}{=}\beta_{\theta,1}M_{\theta}+(1-\beta
_{\theta,1})V$.
\end{longlist}
\end{theorem}
\begin{pf}
From \cite{BondBook}, Section 9, we know that $\Omega_{1}$
corresponds to a $\operatorname{GGC}(1,V)$ subordinator if $q(x)=F_{1/V}(x)$.
Hence, scaling by $\delta$, and using known properties of
variables in $\mathcal{G}_{+}$ concludes the result.
\end{pf}

Note from,\vspace*{-1pt} for instance, \cite{BondBook}, Example 9.2.3, it follows
that for each $\delta>0$, one can choose
$\Omega_{\delta}(1)\stackrel{d}{=}-{\log}(\beta_{\delta,\kappa
-\delta})$
or
$\Omega_{\delta}(1)\stackrel{d}{=}-{\log}(\beta_{\kappa,1+\delta
-\kappa})$
for the beta variables in Example \ref{betaexample}
satisfying (\ref{BArcsine}). Furthermore, among these, the only
choice corresponding to a GGC variable is
$\Omega_{\delta}(1)\stackrel{d}{=}-{\log}(\beta_{\delta,1})$. However,
Theorem \ref{decompositionG} allows us to construct many
quantile clocks based on variables in $\mathcal{G}_{+}$ whose
distributional properties are explicit. We next describe an
interesting property of the variable $ \Sigma_{1}(p)$.
\begin{prop}Suppose that $T_{R}$ is a quantile clock with parameters
$(R,L)$ such that, for an independent variable $Y$,
\[
RY\stackrel{d}{=}U^{1/\delta}\stackrel{d}{=}\beta_{\delta
,1}\stackrel{d}{=}U^{p}
\]
for $\delta=1/p>1$. If $L$ is chosen such that
%
%e6.4 ###
\begin{equation} \label{Lspec2}
\psi_{L}(\omega)=\mathbb{E}\bigl[\psi_{Z}\bigl(\omega
Y{\rme}^{-\Sigma_{1}(p)}\bigr)\bigr],
\end{equation}
where $Z$ is a $\mathcal{U}_{1}$-BDLP, satisfying (\ref{Uvanilla})
for a subordinator $\zeta\in{\mathcal U}_{1}$, then
$T_{R}(t)\stackrel{d}{=}\zeta(t)$ for each fixed $t$.
\end{prop}
\begin{pf}Setting $V=p$, it follows from (\ref{GGCdecomp}), with
$\theta=1$, that
\[
\gamma_{1}\stackrel{d}{=}p\gamma_{1}+\Sigma_{1}(p),
\]
which gives the identity
\[
U^{p}{\rme}^{-\Sigma_{1}(p)}\stackrel{d}{=}U.
\]
Hence, (\ref{Lspec2}) leads to $\mathbb{E}[\psi_{L}(\omega
R)]=\mathbb{E}[\psi_{Z}(\omega U)]$.
\end{pf}

%s6.1 ###
\subsection{Interpreting $\Sigma_{t}(V)$ via diffusions straddling an
exponential time}\label{sec6.1}

Provided that $0<\mathbb{E}[-{\log}(V)]<\infty$, the random
variable $\Delta_{\mathbf{e}}(X)$, with density defined
by (\ref{straddle}) has an interesting interpretation that we now
discuss. This will also give us an opportunity to describe some
more explicit examples of $Q_{R}$. Let $\mathbf{e}/\tilde{X}$
denote an independent $\operatorname{exponential}(1)$ time $\mathbf{e}$ divided
by an independent variable $\tilde{X}$ with distribution
characterized, for bounded measurable functions $H$, by
\[
\mathbb{E}[H(\tilde{X})]=\mathbb{E}[H(X)\log(1+X)]/\varrho.
\]
Now let ${\{\mathcal{R}^{(0,1)}_{s}, s\ge0\}}$ denote a recurrent
linear diffusion starting at $0$ whose inverse local time, in this
case, is a $\operatorname{gamma}(1)$ subordinator. Define for any $t>0$,
\[
g_{t} := \sup \bigl\{s \le t  ;  \mathcal{R}^{(0,1)}_{s}
=0\bigr\},\qquad
 d_{t} := \inf\bigl\{s \ge t,
\mathcal{R}^{(0,1)}_{s}=0\bigr\}.
\]
Then given $\tilde{X}=\lambda$, it follows from  (\ref{straddle})
that for an independent $\operatorname{exponential}(\lambda)$ variable
$\mathbf{e}/\lambda$, the random variable
\[
\Delta_{\mathbf{e}}(\lambda)\stackrel{d}{=}d_{\mathbf{e}/\lambda
}-g_{\mathbf{e}/\lambda}
\]
corresponds to the length of excursion of $\mathcal{R}^{(0,1)}$
above $0$ straddling an\break $\operatorname{exponential}(\lambda)$ time. See, for
instance, \cite{Salminen}, Section 4, for this description for more
general $\mathcal{R}$ as well as \cite{BFRY,Winkel,Pit97}. In
addition, see \cite{BFRY,JRY,JamesBernoulli} for
$\Delta_{\mathbf{e}}(\lambda)$ representation as a variable in
$\mathcal{G}_{+}$. Hence, $\Delta_{\mathbf{e}}(X)$ interprets as
$\Delta_{\mathbf{e}}(\lambda)$ but now for a random time
$\mathbf{e}/\tilde{X}$, with c.d.f. $F_{\mathbf{e}/\tilde{X}}$
satisfying
\[
1-F_{\mathbf{e}/\tilde{X}}(y)=\mathbb{E}[{\rme}^{-Xy}\log(1+X)]/\varrho.
\]

It follows that for $\lambda=(1-p)/p$ and $\varrho=-{\log}(p)$,
that
\[
\Sigma_{t}(p)\stackrel{d}{=}\sum_{k=1}^{N(\varrho
t)}d^{(k)}_{\mathbf{e}/\lambda}-g^{(k)}_{\mathbf{e}/\lambda},\qquad
t>0,
\]
where
$(d^{(k)}_{\mathbf{e}/\lambda},g^{(k)}_{\mathbf{e}/\lambda})$ are
i.i.d. copies of $(d_{\mathbf{e}/\lambda},g_{\mathbf{e}/\lambda})$.

%s6.2 ###
\subsection{Some related examples}\label{sec6.2}

From the results in \cite{BFRY} (see also
\cite{JRY,JamesBernoulli}), we can consider more generally
$\mathcal{R}^{(\alpha,1)}$, in place of $\mathcal{R}^{(0,1)}$,
which, for $0\leq\alpha<1$ is now a process whose inverse local
time is distributed as a generalized gamma subordinator with
L\'{e}vy density specified by $s^{-\alpha-1}{\rme}^{-s}/\Gamma(1-\alpha)$ for $s>0$. Furthermore, for $\lambda=1$,
$\Delta^{(\alpha,1)}_{\mathbf{e}}$, is the generalization of
$\Delta_{\mathbf{e}}(1)\stackrel{d}{=}\Delta^{(0,1)}_{\mathbf{e}}$,
with density
%
%e6.5 ###
\begin{equation}
\label{dena}
\Delta^{(\alpha,1)}_{\mathbf{e}}\stackrel{d}{=}\frac{\alpha
x^{-\alpha-1}{\rme}^{-x}(1-{\rme}^{-x})}{[2^{\alpha}-1]\Gamma(1-\alpha)}\qquad
\mbox{for }x>0.
\end{equation}
Note that the variable $U_{\alpha,\mathbf{e}}\stackrel{d}{=}{\rme}^{-\Delta^{(\alpha,1)}_{\mathbf{e}}}$ has density
%
%e6.6 ###
\begin{equation}
\label{uena} f_{U_{\alpha,\mathbf{e}}}(u)\stackrel{d}{=}\frac
{\alpha
{[-{\log}(u)]}^{-\alpha-1}(1-u)}{[2^{\alpha}-1]\Gamma(1-\alpha
)}\qquad\mbox{for }0<u\leq1.
\end{equation}
In addition, \cite{BFRY} show that
$\Delta^{(\alpha,1)}_{\mathbf{e}}$ is
$\operatorname{GGC}(1-\alpha,\mathbb{D}_{\alpha})$, where $\mathbb{D}_{\alpha}$
satisfies
\[
\mathbb{G}_{\alpha}\stackrel{d}{=}\frac{1}{\mathbb{D}_{\alpha}}-1
\]
with
\[
\log(X_{1-\alpha})=\log(S_{1-\alpha}/S'_{1-\alpha})=\frac{\alpha
}{1-\alpha}\log\bigl(\mathbb{G}_{\alpha}/(1-\mathbb{G}_{\alpha})\bigr).
\]
Furthermore, $\mathbb{G}_{1/2}\stackrel{d}{=}\beta_{1/2,1/2}$,
$\mathbb{G}_{1}\stackrel{d}{=}U$ and
$1/\mathbb{G}_{0}\stackrel{d}{=}1+{\rme}^{\pi\eta}$ for
$\eta$ a standard Cauchy variable. Furthermore
$\gamma_{1-\alpha}U\stackrel{d}{=}\gamma_{1}\beta_{1-\alpha
,1+\alpha}$
is $\operatorname{GGC}(1-\alpha,\mathbb{G}_{\alpha})$. We now look at some
special case of Theorem \ref{decompositionG}.
\begin{prop}\label{straddecomp}Let $0\leq\alpha<1$.
\begin{longlist}[(ii)]
\item[(i)]Then for $V=\mathbb{D}_{\alpha}$ and $X=\mathbb
{G}_{\alpha}$,
%
%e6.7 ###
\begin{eqnarray}
\gamma_{1}&\stackrel{d}{=}&\Delta^{(\alpha,1)}_{\mathbf{e}}+\Sigma
_{1-\alpha}(\mathbb{D}_{\alpha})+\gamma_{\alpha}
\nonumber\\
&\stackrel{d}{=}&\Delta^{(\alpha,1)}_{\mathbf{e}}+\gamma_{\alpha
}M_{\alpha}+\Sigma_{1}(\mathbb{D}_{\alpha})\\
&\stackrel{d}{=}&\Delta^{(\alpha,1)}_{\mathbf{e}}+\Sigma_{1}(\xi
_{1-\alpha}\mathbb{D}_{\alpha}),\nonumber
\end{eqnarray}
where $\Delta^{(\alpha,1)}_{\mathbf{e}}$ has density (\ref{dena}).
$M_{\alpha}\stackrel{d}{=}\beta_{\alpha,1}M_{\alpha}+(1-\beta
_{\alpha,1})\mathbb{D}_{\alpha}$.
When $\alpha=1/2$,
$\gamma_{1/2}M_{1/2}\stackrel{d}{=}\Delta^{(1/2,1)}_{\mathbf{e}}$,
otherwise, $\gamma_{\alpha}M_{\alpha}$ has an explicit density
given in \cite{JamesBernoulli}, Theorem 4.2.
\item[(ii)]For $V=\mathbb{G}_{\alpha}$ and $X=(X_{1-\alpha})^{
{(1-\alpha)}/{\alpha}}$,
%
%e6.8 ###
\begin{eqnarray}
\label{specialG}
\gamma_{1}&\stackrel{d}{=}&\gamma_{1-\alpha}U+\Sigma_{1-\alpha
}(\mathbb{G}_{\alpha})+\gamma_{\alpha}
\nonumber\\
&\stackrel{d}{=}&\gamma_{1-\alpha}U+\gamma_{\alpha}\tilde
{M}_{\alpha}+\Sigma_{1}(\mathbb{G}_{\alpha})\\
&\stackrel{d}{=}&\gamma_{1-\alpha}U+\Sigma_{1}(\xi_{1-\alpha
}\mathbb{G}_{\alpha}),\nonumber
\end{eqnarray}
where
$\tilde{M}_{\alpha}\stackrel{d}{=}\beta_{\alpha,1}\tilde
{M}_{\alpha}+(1-\beta_{\alpha,1})\mathbb{G}_{\alpha}$.
\end{longlist}
\end{prop}
\begin{pf}Note that
$\Delta^{(\alpha,1)}_{\mathbf{e}}+\gamma_{\alpha}M_{\alpha}$ is
$\operatorname{GGC}(1,\mathbb{D}_{\alpha})$ and the second equality above is a
direct consequence of Theorem \ref{decomposition}. Additionally,
it follows that since $\Delta^{(\alpha,1)}_{\mathbf{e}}$ is
$\operatorname{GGC}(1-\alpha,\mathbb{D}_{\alpha})$,
$\gamma_{1-\alpha}\stackrel{d}{=}
\Delta^{(\alpha,1)}_{\mathbf{e}}+\Sigma_{1-\alpha}(\mathbb
{D}_{\alpha})$.
The last equality follows from $\Delta^{(\alpha,1)}_{\mathbf{e}}$
is $\operatorname{GGC}(1,\xi_{1-\alpha}\mathbb{D}_{\alpha})$.
\end{pf}

Recall the\vspace*{-1pt} occupation time variables in
Example \ref{Arcsineexample}, where the variable
$\tilde{\mathbb{O}}_{\alpha,p}$ has density (\ref{denB}), and
additionally $\mathbb{O}_{0}=1/2$ and
$\mathbb{O}_{1}\stackrel{d}{=}\xi_{1/2}$. Then using
Theorems \ref{decomposition} or \ref{decompositionG}, we
obtain the following result.
\begin{prop}Let $0\leq\alpha\leq1$ and $0<p\leq1$.
\begin{longlist}[(ii)]
\item[(i)]Then for $V=\mathbb{O}_{\alpha}$ and $X=X_{\alpha}$,
%
%e6.9 ###
\begin{equation}
\label{specialG2}
\gamma_{1} \stackrel{d}{=} \gamma_{1}\tilde{\mathbb
{O}}_{\alpha,p}+\Sigma_{p}(\mathbb{O}_{\alpha})+\gamma_{1-p}.
\end{equation}
\item[(ii)]When $\alpha=1/2$, $V\stackrel{d}{=}\beta_{1/2,1/2}$ and
$X\stackrel{d}{=}\gamma_{1/2}/\gamma'_{1/2}$,
%
%e6.10 ###
\begin{equation}
\label{specialG3}
\gamma_{1}\stackrel{d}{=}\gamma_{p}\beta_{p+1/2,p+1/2}+\Sigma
_{p}(\beta_{1/2,1/2})+\gamma'_{1-p}.
\end{equation}
\end{longlist}
\end{prop}
\begin{rem} The diffusions above belong to a more general
family, $\mathcal{R}^{(\alpha,b)}$ for $0\leq\alpha<1, b\ge0$
with inverse local time corresponding to a generalized gamma
subordinator with L\'{e}vy density $Cs^{-\alpha-1}{\rme}^{-bs}$.
Hence, the density of variables $\Delta^{(\alpha,b)}_{\mathbf{e}}$
is proportional to $s^{-\alpha-1}{\rme}^{-bs}(1-e^{-s})$. In
particular, for $b=0$, the variable
$\Delta^{(\alpha,0)}_{\mathbf{e}}\stackrel{d}{=}\gamma_{1-\alpha
}/U^{1/\alpha}$
is a $\operatorname{GGC}(1-\alpha,1/\mathbb{G}_{\alpha})$ variable.
See \cite{BFRY,JRY,JamesBernoulli} for more details.
\end{rem}
\begin{rem}\label{remarkexample}
Reference \cite{BFRY}, Theorem 1.4, yields the following
decomposition of $\gamma_{1}$, for $0\leq\alpha\leq1$:
\[
\gamma_{1}\stackrel{d}{=}\gamma_{1}\mathbb{G}_{\alpha}+\gamma
'_{1}\mathbb{G}_{1-\alpha}.
\]
As a special case, with $\alpha=0$ or $1$,
\[
\gamma_{1}\stackrel{d}{=}\gamma_{1}U+\gamma'_{1}/(1+{\rme}^{\pi
\eta}).
\]
Combining this fact with with (\ref{GGCdecomp}), with
$V\stackrel{d}{=}\mathbb{G}_{0}\stackrel{d}{=}1/(1+{\rme}^{\pi
\eta})$, leads to the interesting identity
\[
\Sigma_{1}(\mathbb{G}_{0})\stackrel{d}{=}\gamma_{1}\mathbb{G}_{0}
\]
with L\'{e}vy density given by (\ref{straddle}). This follows since
$\gamma_{1}U$ is a $\operatorname{GGC}(1,\mathbb{G}_{0})$ variable.
\end{rem}

%s7 ###
\section{Composition of quantile clocks}\label{sec7}

We now highlight an important property of the general class of
quantile clocks. Recall from Proposition \ref{Qclock1} that for
each fixed $t$, the marginal distribution of a quantile clock
$T_{R}$, with parameters $(R,L)$, satisfies
\[
T_{R}(t)\stackrel{d}{=}\zeta(t),
\]
where $\zeta$ is a subordinator such that
$\zeta(1)$ has Laplace exponent
\[
\psi_{\zeta}(\omega)=\mathbb{E}[\psi_{L}(\omega
R)]=\psi_{T_{R}(1)}(\omega).
\]
An important operation for L\'{e}vy processes is the composition of
L\'{e}vy processes, in financial applications this is associated
with time changed processes. The fact that quantile clocks behave
marginally like a subordinator allows us to obtain the following
results.

First, we can discuss the composition of two independent quantile
clocks $T_{R_{1}}$, $T_{R_{2}}$, with parameters $(R_{1},L_{1})$
and $(R_{2},L_{2})$, respectively, which can be written as
%
%e7.1 ###
\begin{equation}
\label{compclock}
T_{R_{1}}(T_{R_{2}}(t))=\int_{0}^{T_{R_{2}}(t)}Q_{R_{1}}\biggl(\biggl(1-\frac
{s}{T_{R_{2}}(t)}\biggr)_{+}\biggr)L_{1}(ds)
\end{equation}
for $t\ge0$. The apparently complicated random process appearing
in (\ref{compclock}) is no longer a quantile clock. However,
as the next result shows, its marginals are easy to describe. We
use the notation $\circ$ to denote the composition of functions,
so, for instance, $T_{R_{1}}\circ T_{R_{2}}$ means the operation
in (\ref{compclock}).
\begin{prop} Let\vspace*{-1pt} $T_{R_{i}}$, $i=1,\ldots, k$, denote independent
quantile clocks
such that pointwise $T_{R_{i}}(t)\stackrel{d}{=}\zeta_{i}(t)$ for
independent subordinators with corresponding Laplace exponents
$\psi_{\zeta_{i}}(\omega)$ for $i=1,\ldots, k$. Then the\vspace*{1pt}
composition
$\widehat{T}_{k}:=(\widehat{T}_{k}(t)=T_{R_{1}}\circ\cdots\circ
T_{R_{k}}(t)\dvtx t\ge0)$ is an increasing process such that
for each fixed~$t$,
\[
\widehat{T}_{k}(t)\stackrel{d}{=}\widehat{\zeta}_{k}(t)=\zeta
_{1}\circ\cdots\circ
\zeta_{k}(t),
\]
where $\widehat{\zeta}_{k}$ is a subordinator with Laplace exponent
\[
\psi_{\zeta_{k}}\circ\cdots\circ\psi_{\zeta_{1}}(\omega).
\]
If each $T_{R_{i}}$ is a continuous process, then $\widehat{T}_{k}$
is a continuous process.
\end{prop}
\begin{pf} It suffices to show this for $k=2$. But this is
immediate from Proposition \ref{Qclock1}, since given $T_{R_{2}}$,
\[
\psi_{T_{R_{1}}(T_{R_{2}}(t))}(\omega)=T_{R_{2}}(t)\psi_{\zeta
_{1}}(\omega).
\]
\upqed\end{pf}

Of course one can also compose these clocks with subordinators as
follows. The next result is immediate.
\begin{prop}Let $T_{R}$ denote a quantile clock that satisfies
$T_{R}(t)\stackrel{d}{=}\zeta_{1}(t)$ for some subordinator
$\zeta_{1}$. Furthermore, let $\zeta_{2}$ denote a subordinator
independent of $T_{R}$ and $\zeta_{1}$. Then for each fixed $t$
\[
T_{R}(\zeta_{2}(t))\stackrel{d}{=}\zeta_{1}(\zeta_{2}(t))\quad\mbox{and}
\quad\zeta_{2}(T_{R}(t))\stackrel{d}{=}\zeta_{2}(\zeta_{1}(t)).
\]
\end{prop}

We now illustrate an important special case.
\begin{prop}\label{comptiltedstable}For $0<\alpha<1$ and $0<\beta
<1$, one can use the
specifications in Theorem \ref{SelfDecomposableThm} to construct
independent clocks $T_{R_{1}}$ and $T_{R_{2}}$ such that
marginally $T_{R_{1}}(t)\stackrel{d}{=}\widehat{S}_{\alpha}(t)$ and
$T_{R_{2}}(t)\stackrel{d}{=}\widehat{S}_{\beta}(t)$, where
$\widehat{S}_{\alpha}$ and $\widehat{S}_{\beta}$ are independent with
Laplace exponents $[{(1+\omega)}^{\alpha}-1]$ and
$[{(1+\omega)}^{\beta}-1]$. Then for each fixed $t$
\[
T_{R_{1}}(T_{R_{2}}(t))\stackrel{d}{=}T_{R_{2}}(T_{R_{1}}(t))\stackrel
{d}{=}T_{R_{1}}(\widehat{S}_{\beta}(t))\stackrel{d}{=}\widehat{S}_{\alpha
}(T_{R_{2}}(t))
\stackrel{d}{=}\widehat{S}_{\alpha\beta}(t),
\]
that is, for each fixed $t$ the Laplace exponent is
$t[{(1+\omega)}^{\alpha\beta}-1]$. Additionally, the first two
compositions can be specified such that the resulting processes
are continuous, but the latter compositions always correspond to
processes with jumps.
\end{prop}

%s8 ###
\section{Continuous VG, CGMY, NIG and other price processes}\label{sec8}

Summarizing, we have demonstrated that quantile clocks $T_{R}$ can
either be chosen to have strictly continuous and increasing paths or
can be expressed as $T_{\tilde{R}}(t)+aL(t)$, where $T_{\tilde{R}}$ is
a continuous increasing quantile clock. In general, quantile clocks
have marginals that are equivalent to those of a subordinator, $\zeta$,
for each $t$, that is, $T_{R}(t)\stackrel{d}{=}\zeta(t)$. Moreover, we
have shown that for a large class of quantiles $Q_{R}$ we can choose
$T_{R}$ to have any desired marginal law in $\mathcal{U}_{\delta}$, by
choosing a random variable $Y$ and the subordinator $L$ in a clearly
prescribed fashion. Furthermore, our results in the last section show
that composition operations involving quantile clocks or quantile
clocks with subordinators are marginally equivalent in distribution to
compositions of subordinators. All these properties make them highly
desirable components in pricing models based on time changes. For
example, the processes
\[
\bigl(\widehat{W}_{\mu}(T_{R}(t))\dvtx t\ge0\bigr)\quad\mbox{and}\quad
\bigl(\widehat{W}_{\mu}(T_{R_{1}}(T_{R_{2}}(t)))\dvtx t \ge0\bigr)
\]
can be chosen such that they are processes with continuous
trajectories, but have simple and familiar marginal laws. In
addition, for a subordinator $\tilde{\zeta}$, the processes
\[
\bigl(\widehat{W}_{\mu}(\tilde{\zeta}(T_{R}(t)))\dvtx t \ge0\bigr)
\quad\mbox{and}\quad\bigl(\widehat{W}_{\mu}(T_{R}(\tilde{\zeta}(t)))\dvtx t\ge0\bigr)
\]
have jumps, exhibit volatility clustering, and otherwise may be
chosen to have familiar marginal distributions, in fact the same
marginal, for many choices of~$Q_{R}$. We illustrate these points
through some examples that equate these processes marginally with
some of the most popular L\'{e}vy processes.
\begin{example}[(Continuous variance gamma processes)]\label{VGexample}
As a first example, it follows from III of
Theorem \ref{SelfDecomposableThm}, that if $RY\stackrel{d}{=}U$, then
for each $\delta>0$, a quantile clock $T_{R^{1/\delta}}$ with
parameters $(R^{1/\delta},L_{\delta})$, can be chosen such that
for each fixed $t$,
\[
\widehat{W}_{\mu}(T_{R^{1/\delta}}(t))\stackrel{d}{=}\widehat
{W}_{\mu}(\gamma_{\theta}(t)),
\]
that is, it has marginal distributions equivalent to the log
price of a variance gamma (VG) process \cite{MadanVG}, not
depending on $\delta$, if for each $\delta>0$, the subordinator
(depending on $\delta$)
\[
L_{\delta}(s)\stackrel{d}{=}\zeta_{\delta,Y^{1/\delta}}(s)+\sum
_{k=1}^{N(\theta
s/\delta)}\gamma^{(k)}_{1}Y^{1/\delta}_{k},\qquad s\ge0,
\]
where $\zeta_{\delta,Y^{1/\delta}}$ is a $\operatorname{GGC}(\theta,
Y^{1/\delta})$ subordinator.
\end{example}

We next show how to obtain price processes whose marginal laws are
equivalent to a Carr--Geman--Madan--Yor (CGMY) process \cite{CGMY}
but otherwise possesses continuous sample paths.
\begin{example}[(Continuous CGMY processes)]\label{CGMYexample}
For this example, we follow the exposition in \cite{MadanYorCGMY}.
Let
\[
A=\frac{G-M}{2}\quad\mbox{and}\quad B=\frac{G+M}{2},
\]
then the L\'{e}vy density of the log prices of a CGMY process, say
$\chi_{\mathrm{CGMY}}$, is given by
%
%e8.1 ###
\begin{equation}
\label{levyCGMY}\qquad
\rho_{\chi_{\mathrm{CGMY}}(1)}(x)=\frac{\Gamma(\alpha)\Gamma(1-\alpha
)}{\Gamma(1+\alpha)}{\rme}^{Ax-B|x|}|x|^{-d-1}\qquad
\mbox{for }-\infty<x<\infty
\end{equation}
and $0<d=2\alpha<2$. Madan and Yor \cite{MadanYorCGMY}, show that
the log price of a CGMY process has an explicit representation in
terms of a time changed brownian motion,
$\chi_{\mathrm{CGMY}}(t):=\widehat{W}_{A}(\zeta(t))$, where $\zeta$ is a
subordinator with L\'{e}vy density
%
%e8.2 ###
\begin{eqnarray}
\label{levyCGMYsub}
\rho_{\zeta}(s)&=&\frac{2^{\alpha}\Gamma(\alpha)}{\Gamma(2\alpha)}{\rme}^{{(A^{2}-B^{2})s}/{2}}s^{-\alpha-1}
\mathbb{E}\bigl[{\rme}^{-s({B^2}/{2})({\gamma_{\alpha}}/{\gamma_{1/2}})}\bigr]\nonumber\\[-8pt]\\[-8pt]
&=&\frac{2^{\alpha}\Gamma(\alpha)}{\Gamma(2\alpha)}s^{-\alpha
-1}\mathbb{E}[{\rme}^{-sV}]\nonumber
\end{eqnarray}
for
%
%e8.3 ###
\begin{equation}
\label{VCGMY}
V\stackrel{d}{=}\biggl(4MG+B^2\frac{\gamma_{\alpha}}{\gamma_{1/2}}\biggr)\Big/2.
\end{equation}
It is evident from (\ref{levyCGMYsub}) that $\zeta\in
\mathcal{G}$. We now give the specifications for a quantile clock
to have marginals with L\'{e}vy density (\ref{levyCGMYsub}) hence
inducing price processes that have the marginal distribution of a
CGMY process. Notice that
%
%e8.4 ###
\begin{eqnarray}
\label{BDLPCGMY}
-s\rho'_{\zeta}(s)&=&\frac{2^{\alpha}\Gamma(\alpha)}{\Gamma
(2\alpha)}s^{-\alpha-1}
\mathbb{E}\bigl[\bigl((1+\alpha)+sV\bigr){\rme}^{-sV}\bigr]\nonumber\\[-8pt]\\[-8pt]
&=&(1+\alpha)\rho_{\zeta}(s)+\frac{2^{\alpha}\Gamma(\alpha
)}{\Gamma(2\alpha)}
s^{-\alpha}\mathbb{E}[V{\rme}^{-sV}].\nonumber
\end{eqnarray}
Hence, $ \psi_{\zeta}(\omega)=\mathbb{E}[\psi_{Z}(\omega
U^{1/\delta})]$, for
%
%e8.5 ###
\begin{equation}
\label{ZCGMY}
\rho_{Z}(s)=(1+\alpha/\delta)\rho_{\zeta}(s)+\frac{2^{\alpha
}\Gamma(\alpha)}{\delta\Gamma(2\alpha)}
s^{-\alpha}\mathbb{E}[V{\rme}^{-sV}].
\end{equation}
\begin{prop}\label{CGMYspec}Suppose that $T_{R}$ is a quantile clock
with parameters
$(R,L)$, such that there exists a variable $Y$ satisfying
$RY\stackrel{d}{=}U^{1/\delta}$ for some $\delta>0$. Then for each
fixed $t$,
\[
\widehat{W}_{A}(T_{R}(t))\stackrel{d}{=}\chi_{\mathrm{CGMY}}(t)
\]
specified by (\ref{levyCGMY}) if the subordinator $L$ is chosen
such that
%
%e8.6 ###
\begin{eqnarray}
\label{CGMYL}
\rho_{L}(s)&=&(1+\alpha/\delta)\frac{c_{\alpha}2^{\alpha}\Gamma
(\alpha)}{\Gamma(2\alpha)}s^{-\alpha-1}\mathbb{E}[{\rme}^{-sV/Y_{\alpha}}]\nonumber\\[-8pt]\\[-8pt]
&&{}+\frac{c_{\alpha}2^{\alpha}\Gamma(\alpha)}{\delta\Gamma
(2\alpha)}s^{-\alpha}\mathbb{E}[(V/Y_{\alpha}){\rme}^{-sV/Y_{\alpha}}],\nonumber
\end{eqnarray}
where $V$ is defined by (\ref{VCGMY}),
$c_{\alpha}=\mathbb{E}[Y^{\alpha}]$ and $Y_{\alpha}$ is the random
variable whose distribution is proportional to
$y^{\alpha}F_{Y}(dy)$. When $Y=1$, $L:=Z$
satisfying (\ref{ZCGMY}). Note also that $\mathbb{E}[V^{\alpha}]$
is finite only if $\alpha<1/2$. Hence,
\[
s^{-\alpha}\mathbb{E}[(V/Y_{\alpha}){\rme}^{-sV/Y_{\alpha}}]
\]
is the L\'{e}vy density of a compound Poisson process only in the
case where $\alpha<1/2$.
\end{prop}
\begin{pf}The result is a special case of Theorem \ref
{SelfDecomposableThm} and the specifications we derived above. In
particular, (\ref{ZCGMY}).
\end{pf}

As a specific example with continuous paths, consider again the
Kumaraswamy quantile clock with
\[
T_{K_{p,b}}(t)=\int_{0}^{t}\bigl[1-\bigl(1-(1-s/t)_{+}\bigr)^{1/p}\bigr]^{1/b}L(ds)
\]
for
$R\stackrel{d}{=}K_{p,b}\stackrel{d}{=}(1-U^{1/p})^{1/b}\stackrel
{d}{=}\beta^{1/b}_{1,p}$.
Then for each fixed $t$,
\[
\widehat{W}_{A}(T_{K_{p,b}}(t))\stackrel{d}{=}\chi_{\mathrm{CGMY}}(t),
\]
if $L$ is selected according to (\ref{CGMYL}) with
$Y\stackrel{d}{=}\beta^{1/b}_{p,1-p}$ and $\delta=bp$.

If we consider the arcsine clock using (\ref{TrueArcsine}) with
\[
T_{\beta^{1/b}_{1/2,1/2}}(t):=\int_{0}^{t}\sin^{2/b}\biggl(\frac{\pi
}{2}\biggl(1-\frac{s}{t}\biggr)_{+}\biggr)L(ds)
\]
then for each fixed $t$,
\[
\widehat{W}_{A}(T_{\beta^{1/b}_{1/2,1/2}}(t))\stackrel{d}{=}\chi_{\mathrm{CGMY}}(t),
\]
if $L$ is selected according to (\ref{CGMYL}) with
$Y\stackrel{d}{=}{(1-U^{2})}^{1/b}$ and $\delta=b/2$.

The next two cases are from Section \ref{sec6.2} involving quantile
functions that can be evaluated numerically. If we consider the
clock $T_{U^{1/\delta}_{\alpha,\mathbf{e}}}$ based on the
variable $U_{\alpha,\mathbf{e}}$ with density (\ref{uena}), it
follows from Proposition \ref{straddecomp} that
\[
\widehat{W}_{A}(T_{U^{1/\delta}_{\alpha,\mathbf{e}}}(t))\stackrel
{d}{=}\chi_{\mathrm{CGMY}}(t),
\]
if $L$ is selected
according to (\ref{CGMYL}) with
\[
Y\stackrel{d}{=}{\rme}^{-[\Sigma_{1-\alpha}(\mathbb{D}_{\alpha})+\gamma_{\alpha
}]/\delta}.
\]
If we consider the variables in Remark \ref{remarkexample} then this
leads to a quantile clock based on the variable $\gamma_{1}\mathbb
{G}_{\alpha}$. Hence,
\[
\widehat{W}_{A}(T_{{\rme}^{-\gamma_{1}\mathbb{G}_{\alpha}/\delta
}}(t))\stackrel{d}{=}\chi_{\mathrm{CGMY}}(t),
\]
if $L$ is selected
according to (\ref{CGMYL}) with
\[
Y\stackrel{d}{=}{\rme}^{-\gamma_{1}\mathbb{G}_{1-\alpha}/\delta}.
\]
Other examples
using (\ref{specialG2}) and (\ref{specialG3}) are based on the
pairs
\[
\bigl({\rme}^{- \gamma_{1}\tilde{\mathbb{O}}_{\alpha,p}}, {\rme}^{-[\Sigma_{p}(\mathbb{O}_{\alpha})+\gamma_{1-p}]}\bigr)
\quad\mbox{and}\quad \bigl({\rme}^{-\gamma_{p}\beta_{p+1/2,p+1/2}}, {\rme}^{-[\Sigma_{p}(\beta_{1/2,1/2})+\gamma'_{1-p}]}\bigr).
\]

Finally, if instead one uses the clock
\[
T_{\tilde{U}^{1/\delta}_{p}}(t):=\int_{0}^{t}{\biggl[(1-p)+p\biggl(1-\frac
{s}{t}\biggr)_{+}\biggr]}^{1/\delta}L(ds),
\]
then for each fixed $t$,
\[
\widehat{W}_{A}(T_{\tilde{U}^{1/\delta}_{p}}(t))\stackrel{d}{=}\chi
_{\mathrm{CGMY}}(t),
\]
if $L$ is selected according to (\ref{CGMYL}) with $Y={\rme}^{-X_{p}[-{\log}(1-p)]/\delta}$ where again $X_{p}$ is
$\operatorname{geometric}(p)$. Hence, for each $0<p<1$, the resulting process has
CGMY marginals, exhibits volatility clustering, but also has
jumps. If $p=1$, then the process is continuous and the quantile
clock coincides with the Holmgren--Liouville clock discussed in
\cite{Bender}.
\end{example}

It is evident that the specifications (\ref{CGMYL}) for $L$
appearing in Proposition \ref{CGMYspec} can be modified such that
$T_{R}(t)$ has marginals equivalent to a subordinator with L\'{e}vy
density
\[
\rho_{\zeta}(s)=Cs^{-\alpha-1}\mathbb{E}[{\rme}^{-sV}]
\]
for some positive constant $C$, where $V$ is a much more general
random variable. That is, $L$ is specified by
%
%e8.7 ###
\begin{eqnarray}
\label{CGMYLgen}
\rho_{L}(s)&=&(1+\alpha/\delta)c_{\alpha}Cs^{-\alpha-1}\mathbb
{E}[{\rme}^{-sV/Y_{\alpha}}]\nonumber\\[-8pt]\\[-8pt]
&&{}+(c_{\alpha}/\delta)Cs^{-\alpha}\mathbb{E}[(V/Y_{\alpha}){\rme}^{-sV/Y_{\alpha}}].\nonumber
\end{eqnarray}
As a specific example, we next look at the case corresponding to
NIG and related processes.
\begin{example}[(Processes with NIG and related marginals)]\label
{NIGexample}For this example,
let $\widehat{S}_{\alpha}(t)$ denote any subordinator with L\'{e}vy
density
\[
\rho_{\alpha}(s)=\frac{\alpha}{\Gamma(1-\alpha)}s^{-\alpha
-1}{\rme}^{-s}
\]
and define the L\'{e}vy process on $\mathbb{R}$ by
$\chi_{\alpha}(t):=\widehat{W}_{\mu}(\widehat{S}_{\alpha}(t))$. It follows
that
\[
\chi_{1/2}(t):=\chi_{\mathrm{NIG}}(t)
\]
is a normal inverse Gaussian (NIG) process \cite{BN98}. If
$T^{(\alpha)}_{R}$ denotes a quantile clock such that
$RY=U^{1/\delta}$, and $L$ is specified according
to (\ref{CGMYLgen}), specifically using the L\'{e}vy density
$\rho_{\alpha}$, with $V=1$, then
$\widehat{W}_{\mu}(T^{(\alpha)}_{R}(t))\stackrel{d}{=}\chi_{\alpha}(t)$.
In particular setting $\alpha=1/2$, it follows that for each $t$,
\[
\widehat{W}_{\mu}\bigl(T^{(1/2)}_{R}(t)\bigr)\stackrel{d}{=}\chi_{\mathrm{NIG}}(t).
\]
Thus, yielding processes with continuous trajectories but NIG
marginals. In addition, choosing $\alpha$ and $0<\beta\leq1$ such
that $\alpha\beta=1/2$, it follows from
Proposition~\ref{comptiltedstable} that
\[
\chi_{\beta}\bigl(T^{(\alpha)}(t)\bigr)\stackrel{d}{=}\widehat{W}_{\mu
}\bigl(T^{(\alpha)}(\widehat{S}_{\beta}(t))\bigr)\stackrel{d}{=}\chi_{\mathrm{NIG}}(t),
\]
corresponding to processes with jumps, dependent increments and
NIG marginal distributions. Note that when $Y=1$, corresponding to
the quantile clock $T_{U^{1/\delta}}$, then similar to
(\ref{ZCGMY}),
%
%e8.8 ###
\begin{equation}
\label{plainNIG}
L(t)=\widehat{S}_{\alpha}\bigl((1+\alpha/\delta)t\bigr)+\sum_{k=1}^{N(\alpha
s/\delta)}\gamma^{(k)}_{1-\alpha}\qquad\mbox{for }t>0.
\end{equation}
\end{example}

%s9 ###
\section{Choosing laws for the short memory kernel}\label{sec9}

We now apply our results for quantile clocks to a convoluted
subordinator that \cite{Bender} refer to as a short memory
kernel. We note that this convoluted subordinator is not a
quantile clock.
\begin{theorem}Let $\zeta$ denote a subordinator with self-decomposable
laws such that the quantile clock with parameters $(U,Z)$, that is,
$T_{U}$ has marginals $T_{U}(t)\stackrel{d}{=}\zeta(t)\in
\mathcal{L}$. This is achieved by setting the L\'{e}vy density of
$Z$ to be
%
%e9.1 ###
\begin{equation}
\label{Zspecshort}
\rho_{Z}(x)=-x\rho'_{\zeta}(x)=\rho_{\zeta}(x)+\rho_{\vartheta}(x),
\end{equation}
where $\vartheta$ is the OU--BDLP of
$v(0)\stackrel{d}{=}T_{U}(1)\stackrel{d}{=}\zeta(1)$. Then for $Z$
satisfying (\ref{Zspecshort}), the short memory convoluted
subordinator constructed as
\[
\tilde{T}_{\varepsilon}(t)=\int_{0}^{t}\min\biggl(1,\frac
{(t-s)_{+}}{\varepsilon}\biggr)Z(ds)
\]
has the following distributional properties:
\begin{longlist}[(ii)]
\item[(i)]For each fixed $t$, the Laplace exponent of the r.v.
$\tilde{T}_{\varepsilon}(t)$, is given by
\[
\psi_{\tilde{T}_{\varepsilon}(t)}(\omega)=\cases{
t\psi_{\zeta}\biggl(\omega\dfrac{t}{\varepsilon}\biggr), &\quad $t\leq
\varepsilon$,\vspace*{2pt}\cr
t\psi_{\zeta}(\omega)+(t-\varepsilon)\omega\psi'_{\zeta}(\omega
), &\quad $t>\varepsilon$.}
\]
\item[(ii)]For each fixed t, the marginal distribution of $\tilde
{T}_{\varepsilon}(t)$ is
given by
\[
\tilde{T}_{\varepsilon}(t)\stackrel{d}{=}
\cases{
\dfrac{t}{\varepsilon}\zeta(t), &\quad $t\leq\varepsilon$,\vspace*{2pt}\cr
\zeta(t)+\vartheta(t-\varepsilon), &\quad $t>\varepsilon$.}
\]
\end{longlist}
\end{theorem}
\begin{pf}First, notice that in general, for each fixed $t$, the
L\'{e}vy exponent of the random variable
$\tilde{T}_{\varepsilon}(t)$ is given by
%
%e9.2 ###
\begin{equation}
\label{genmemory}
\psi_{\tilde{T}_{\varepsilon}(t)}(\omega)=\int_{0}^{t}\psi
_{Z}\biggl(\omega\min\biggl(1,\frac{(t-s)_{+}}{\varepsilon}\biggr)\biggr)\,ds.
\end{equation}
So for $t\leq\varepsilon$, (\ref{genmemory}) can be expressed as
\[
\int_{0}^{t}\psi_{Z}\biggl(\omega\min\biggl(1,\frac{(t-s)_{+}}{\varepsilon
}\biggr)\biggr)\,ds=t\mathbb{E}\biggl[\psi_{Z}\biggl(\omega\frac{t}{\varepsilon}U\biggr)\biggr]
=t\psi_{\zeta}\biggl(\omega\frac{t}{\varepsilon}\biggr),
\]
where the last equality follows from (\ref{Zspecshort}). For
$t>\varepsilon$, split the interval $[0,t]$ into $[0,t-\varepsilon]$
and $(t-\varepsilon, t]$ then (\ref{genmemory}) becomes
%
%e9.3 ###
\begin{equation}
\label{genI} \varepsilon\mathbb{E}[\psi_{Z}(\omega
U)]+(t-\varepsilon)\psi_{Z}(\omega).
\end{equation}
Now use (\ref{Zspecshort}) to show that (\ref{genI}) is equal to
\[
\varepsilon\psi_{\zeta}(\omega)+(t-\varepsilon)[\psi_{\zeta
}(\omega)+\psi_{\vartheta}(\omega)]
\]
yielding the result.
\end{pf}

Our result now allows one to choose more convenient laws for
$\tilde{T}_{\epsilon}$ which allows one to easily apply the
option pricing formula of \cite{Bender}, as displayed in
Theorem \ref{tina}, either by exact simulation or FFT methods. We
illustrate this in the next example.
\begin{example}[(Short memory convoluted subordinator with NIG related
marginals)]
First, it is interesting to recall from \cite{BNS2001} that the
OU--BDLP, $\vartheta$, leading to
$v(0)\stackrel{d}{=}\widehat{S}_{\alpha}(1)$, as specified in
Example \ref{NIGexample}, has L\'{e}vy density
\[
\rho_{\vartheta}(s)=\frac{\alpha}{\Gamma(1-\alpha)}s^{-\alpha
-1}[\alpha+s]{\rme}^{-s}.
\]
Hence,
\[
\vartheta(s)\stackrel{d}{=}\widehat{S}_{\alpha}(\alpha
s)+\sum_{k=1}^{N(\alpha s)}\gamma^{(k)}_{1-\alpha}\qquad\mbox{for
}s>0.
\]
Now using (\ref{plainNIG}) with $\delta=1$, setting the
subordinator
\[
Z(s)=\widehat{S}_{\alpha}\bigl((1+\alpha)s\bigr)+\sum_{k=1}^{N(\alpha
s)}\gamma^{(k)}_{1-\alpha}\qquad\mbox{for }s>0
\]
leads to the following marginal behavior of the corresponding
short-memory model, for each fixed $t$,
\[
\tilde{T}_{\varepsilon}(t)\stackrel{d}{=}
\cases{\displaystyle \frac{t}{\varepsilon}\widehat{S}_{\alpha}(t), &\quad $t\leq
\varepsilon$,\vspace*{2pt}\cr
\displaystyle \widehat{S}_{\alpha}\bigl(t+\alpha(t-\epsilon)\bigr)+\sum_{k=1}^{N(\alpha
(t-\varepsilon))}\gamma^{(k)}_{1-\alpha}, &\quad $t>\varepsilon$.}
\]
Note however that $\tilde{T}_{\varepsilon}$ is a
continuous process since $Z$ is an infinite activity process that
satisfies the conditions in \cite{Bender}. This leads to price
processes (\ref{priceCS}) with continuous trajectories that have
the following marginal behavior. $ \widehat{W}_{-1/2}(\sigma^{2}
\tilde{T}_{\varepsilon}(t))$ is equivalent in distribution
to
\[
\cases{\displaystyle \widehat{W}_{-1/2}\biggl(\sigma^{2}\frac{t}{\varepsilon}\widehat{S}_{\alpha
}(t)\biggr), &\quad $t\leq\varepsilon$,\vspace*{2pt}\cr
\displaystyle \widehat{W}_{-1/2}\bigl(\sigma^{2}\widehat{S}_{\alpha}\bigl(t+\alpha(t-\epsilon
)\bigr)\bigr)+\widehat{W'}_{-1/2}\Biggl(
\sigma^{2}\sum_{k=1}^{N(\alpha(t-\varepsilon))}\gamma
^{(k)}_{1-\alpha}\Biggr), &\quad $t>\varepsilon$.}
\]
If one sets $\alpha=1/2$, then this reduces to
\[
\cases{\displaystyle \widehat{W}_{-1/2}\biggl(\sigma^{2}\frac{t}{\varepsilon}\widehat
{S}_{1/2}(t)\biggr), &\quad $t\leq\varepsilon$,\vspace*{2pt}\cr
\displaystyle \widehat{W}_{-1/2}\bigl(\sigma^{2}\widehat{S}_{1/2}\bigl((3t-\varepsilon
)/2\bigr)\bigr)+\widehat{W'}_{-1/2}\Biggl(
\sigma^{2}\sum_{k=1}^{N((t-\varepsilon)/2)}\gamma^{(k)}_{1/2}\Biggr), &\quad
$t>\varepsilon$.}
\]
Hence, in this case, for each fixed $t\leq
\varepsilon$, the marginal distribution of the price
process (\ref{priceCS}) follows an NIG distribution, with scale
parameters depending on $t$.
\end{example}

\section*{Acknowledgment}
We would like to thank Professor Marc Yor who informed us about the
decompositions of $\gamma_{1}$ in Examples \ref{exam5.3} and \ref{split}, which also
influenced the development in Section \ref{sec6}.

% imsref loaded by lrinkeviciute, 2011-02-25 08:49:10
% imsref loaded by lrinkeviciute, 2011-02-25 09:03:38

%
\printaddresses

\end{document}